\documentclass{article}

 \PassOptionsToPackage{numbers, compress}{natbib}
\usepackage[preprint]{neurips_2022}
\usepackage[utf8]{inputenc} % allow utf-8 input
\usepackage[T1]{fontenc}    % use 8-bit T1 fonts
\usepackage{hyperref}       % hyperlinks
\usepackage{url}            % simple URL typesetting
\usepackage{booktabs}       % professional-quality tables
\usepackage{amsfonts}       % blackboard math symbols
\usepackage{nicefrac}       % compact symbols for 1/2, etc.
\usepackage{microtype}      % microtypography
\usepackage{xcolor}         % colors
\usepackage{longtable}
\usepackage{booktabs}
\usepackage{graphicx}
\usepackage{extarrows}
\usepackage{subfigure}
\usepackage{textcomp}
\usepackage{booktabs}
\usepackage{multirow}
\usepackage{wrapfig}
\usepackage{graphicx}
\usepackage{enumitem}
\usepackage{soul}
\usepackage{amsmath,amscd,amsbsy,amssymb,amsfonts,latexsym,url,bm,amsthm}
\def\BibTeX{{\rm B\kern-.05em{\sc i\kern-.025em b}\kern-.08em
		T\kern-.1667em\lower.7ex\hbox{E}\kern-.125emX}}
\usepackage{algorithm}
\usepackage{algorithmic}
\usepackage[algo2e,ruled,vlined,boxed,commentsnumbered]{algorithm2e}
\usepackage{threeparttable}

\def\algcomment#1{\textcolor[rgb]{0,0.6,0}{\# #1}}

\newtheorem{definition}{Definition}

\title{A General Framework for Evaluating Robustness of Combinatorial Optimization Solvers on Graphs}

\author{%
  Han Lu\footnotemark[1]~\footnotemark[2] \\
 \And 
 Zenan Li\footnotemark[1]~\footnotemark[2]\\
 \And 
 Runzhong Wang \footnotemark[2]\\
 \And 
 Qibing Ren \footnotemark[2]\\
 \And 
 Junchi Yan \footnotemark[2]\\
 \And
 Xiaokang Yang \footnotemark[2]\\
}
\begin{document}

\maketitle
\renewcommand{\thefootnote}{\fnsymbol{footnote}} %将脚注符号设置为fnsymbol类型，即特殊符号表示
\footnotetext[1]{These authors contributed equally to this work.}
\footnotetext[2]{Department of Computer Science and Engineering, and MoE Key Lab of Artificial Intelligence, AI Institute,
Shanghai Jiao Tong University. Correspondence to: Junchi Yan
<yanjunchi@sjtu.edu.cn>.}

\begin{abstract}
  Solving combinatorial optimization (CO) on graphs is among the fundamental tasks for upper-stream applications in data mining, machine learning and operations research. Despite the inherent NP-hard challenge for CO, heuristics, branch-and-bound, learning-based solvers are developed to tackle CO problems as accurately as possible given limited time budgets. However, a practical metric for the sensitivity of CO solvers remains largely unexplored. Existing theoretical metrics require the optimal solution which is infeasible, and the gradient-based adversarial attack metric from deep learning is not compatible with non-learning solvers that are usually non-differentiable. In this paper, we develop the first practically feasible robustness metric for general combinatorial optimization solvers. We develop a no worse optimal cost guarantee thus do not require optimal solutions, and we tackle the non-differentiable challenge by resorting to black-box adversarial attack methods. Extensive experiments are conducted on 14 unique combinations of solvers and CO problems, and we demonstrate that the performance of state-of-the-art solvers like Gurobi can degenerate by over 20\% under the given time limit bound on the hard instances discovered by our robustness metric, raising concerns about the robustness of combinatorial optimization solvers. 
\end{abstract}

\section{Introduction}\vspace{-5pt}
% KDD paper related
 The combinatorial optimization (CO) problems on graphs are widely studied due to their important applications including aligning cross-modality labels \cite{lyu2020partial}, discovering vital seed users in social networks \cite{maxcovering}, tackling large-scale knapsack problems~\cite{zhang2020solving} and scheduling jobs \cite{mao2019learning} in data centers, etc. However, CO problems are non-trivial to solve due to the NP-hard challenge, whereby the optimal solution can be nearly infeasible to achieve for even medium-sized problems. Existing approaches to practically tackle CO include heuristic methods~\citep{van1987simulated, whitley1994genetic}, powerful branch-and-bound solvers~\citep{llc2020gurobi,scip,cbc} and recently developed learning-based models~\citep{khalil2017learning, mao2019learning, kwon2021matrix}. 

%In the last century, a variety of heuristic methods~\citep{van1987simulated, whitley1994genetic} are proposed to tackle these standing and often NP-hard problems. Driven by the recent development of deep learning and reinforcement learning, many learning-based methods~\citep{khalil2017learning, mao2019learning, kwon2021matrix} are also developed in this area, which show promising potential for their cost-efficiency.

Despite the success of solvers in various combinatorial tasks, little attention has been paid to the vulnerability and robustness of combinatorial solvers. Within the scope of the solvers and problems studied in this paper, our study shows that the performance of the solver may degenerate a lot given certain data distributions that should lead to the same or better solutions compared to the original distribution assuming the solver works robustly. We also validate in experiments that such a performance degradation is neither caused by the inherent discrete nature of CO. Such a discovery raises our concerns about the stability and robustness of combinatorial solvers, which is also aware by \cite{varma2021average,geisler2021generalization}. However, \cite{varma2021average} focuses on theoretical analysis and requires the optimal solution, which is infeasible to reach in practice. \cite{geisler2021generalization} applies the attack-and-defense study in deep learning into learning-based combinatorial optimization models, requiring a differentiable neural network that is infeasible for general solvers that can often be undifferentiable. Our paper also differs from the namely ``robust optimization'' \cite{buchheim2018robust}, where the expected objective score is optimized given known data distribution.
We summarize the challenges and our initiatives of evaluating the robustness of existing solvers as follows: \looseness=-1

\textbf{1) Robustness metric without optimal solutions}: The underlying NP-hard challenge of combinatorial optimization prohibits us from obtaining the optimal solutions. However, the robustness metric proposed by \cite{varma2021average} requires optimal solutions, making this metric infeasible in practice. In this paper, we propose a problem modification based robustness evaluation method, whereby the solvers' performances are evaluated on problems that can guarantee \emph{no worse optimal costs} compared to the original problem, without requiring the optimal solutions.
    
    %Whether the solvers can perform as stably and robustly in real applications as they perform in standard  benchmarks is a common concern \cite{varma2021average,geisler2021generalization}. But how to define the metric of robustness remains a dilemma. NP-hard problems are hard to . A common method is to evaluate solvers using their performance in standard benchmarks. However, the labels for the robustness testing instances are hard to obtain.

\textbf{2) Robustness metric of general (non-differentiable) solvers}: 
    Despite the recently developed deep-learning solvers, most existing solvers are non-differentiable due to the discrete combinatorial nature. Thus, the gradient-based robustness metric by \cite{geisler2021generalization} has the limitation to generalize to general solvers. In this paper, we develop a reinforcement learning (RL) based attacker by regarding arbitrary types of solvers as black-boxes, and the RL agent is trained by policy gradient without requiring the solver to be differentiable. Our framework also owns the flexibility where the agent can be replaced by other search schemes e.g.\ simulated annealing.
    %Some solvers are packaged as a tool and the structure cannot be detected. For example, the commercial solver Gurobi only provides the interfaces while the testers cannot obtain the detailed inner implementation.

% regardless of whether they are learning based or not.  
%A line of relevant works aims at handling combinatorial optimization under uncertainty~\citep{buchheim2018robust} instead of exploring the deterministic problems with slight modification. \cite{varma2021average} initiates a systematic study of average sensitivity on the solution, which is Hamming distance between two outputs by moving an edge. However, these analysis focuses on the deterministic heuristic and often require the ground truth optimal value as the evaluation metric, which is hard to calculate. Meanwhile, the objective of the solution attracts more attention in CO. The recent work \cite{geisler2021generalization} conducts a meaningful attempt on exploring the adversarial robustness on neural combinatorial solvers. The work requires white-box learnable neural solver  and the target CO problems should have specific property to generate the attack instances, which can be calculated the optimal value directly. These shortcomings limit the generalization. 

\begin{table*}[tb!]
    \centering%~\citeyearpar{madry2018towards}
    \caption{Comparing our framework (ROCO) with FGSM~\cite{goodfellow2014explaining} and RL-S2V~\cite{dai2018adversarial}. $\epsilon$-perturb means the change of one pixel should be bounded in $\epsilon$. B-hop neighbourhood means the new attack edges can only connect two nodes with distance less than $B$. }
    \resizebox{\textwidth}{!}
    {
    \begin{tabular}{r|cccccc}
        \hline
        % Table generated by Excel2LaTeX from sheet 'MGM-willow'
        \textbf{Method} & \textbf{Data} & \textbf{Task} & \textbf{Attack target} & \textbf{Attack budget}& \textbf{Attack principle} &  \\
        \hline
        FGSM~\cite{goodfellow2014explaining} & image  & classification  & pixels  & $\epsilon$-Perturb & invisible change    \\
        RL-S2V~\cite{dai2018adversarial} & graph  & classification  & edges (connectivity)& \# edge  & B-hop neighbour   \\
        ROCO (ours) & CO instance  & CO solution  & edges (constraints/cost)  & \# edge & no worse optimum   \\
        \hline
    \end{tabular}%
    }
    \label{tab:different_attack_comparison}
    \vspace{-25pt}
\end{table*}

% To this end, we present \textbf{\underline R}\textbf{\underline o}bust \textbf{\underline C}ombnaotorial \textbf{\underline O}ptimization (ROCO), a general framework to evaluate the robustness of combinatorial optimization solvers. To our best knowledge, this is the first robustness metric for general combinatorial optimization solvers without requiring the optimal solution, making it feasible for us to conduct an extensive evaluation of the robustness of 14 unique combinations of different solvers and problems.

Being aware that many CO problems can be essentially formulated as a graph problem and there are well-developed graph-based learning models~\cite{scarselli2008graph, kipf2016semi, velickovic2017graph}, the scope of this paper is restricted within combinatorial problems on graphs following~\citep{khalil2017learning,wang2021bi}. To develop a problem modification method that achieves \emph{no worse optimal costs}, we propose to modify the graph structure with a problem-dependent strategy exploiting the underlying problem natures. Two types of strategies are developed concerning the studied four combinatorial problems: 1) loosening the constraints such that the feasible space can be enlarged while the original optimal solution is preserved; 2) lowering the cost of a partial problem such that the objective value of the original optimal solution can become better. Such strategies are feasible when certain modifications are performed on the edges of the input graph, and tackle the challenge that we have no access to the optimal solutions. It is also ensured that the generated data distributes closely to the original data by restricting the number of edges modified. The graph modification steps are performed by our proposed attackers that generate problem instances where the solvers' performances degenerate, and we ensure the \emph{no worse optimal costs} guarantee by restricting the action space of the attackers.

%It is worth noting that many CO problems can be essentially formulated as a graph problem~\citep{khalil2017learning,bengio2020machine}, hence it is attractive and natural to modify the problem instance by modifying the graph structure, to generate more test cases for solvers. Therefore we propose the \emph{soft-evaluation} method to modify the instances from the benchmarks and evaluate the robustness where we ensure the modification will generate a better or at least the same optimal solution. The generation process limits the number of modified edges (the constrains in CO) to maintain the similarity between new instances and the original benchmark. The \emph{soft-evaluation} effectively alleviates \emph{No ground truth label} challenge and the reasonable \emph{Robustness metric} derives from it as calculate the difference of the objective function.

Our framework can be seen as the attack on the solvers and how much the worst cases near the clean instances will do harm to the solver are considered as the robustness. Tab.~\ref{tab:different_attack_comparison} compares our framework to classical works on adversarial attacks for images and graphs. 
\textbf{The main contributions of this paper are summarized as follows:}

\textbf{1)} We propose \textbf{\underline R}\textbf{\underline o}bust \textbf{\underline C}ombnaotorial \textbf{\underline O}ptimization (ROCO),  the first general framework (as shown in Fig.~\ref{fig:roco}) to evaluate the robustness of combinatorial solvers on graphs without requiring optimal value or differentiable property.

\textbf{2)} We design a novel robustness evaluation method with \emph{no worse optimal costs} guarantee to eliminate the urgent requirements for the optimal solution of NP-hard CO problems. We develop a reinforcement learning (RL) based attacker combined with other search-based attackers to regard the solvers as black-boxes, making the non-differentiable ones available.

\textbf{3)} Case studies are conducted on four common combinatorial tasks: Directed Acyclic Graph Scheduling, Asymmetric Traveling Salesman Problem, Maximum Coverage, and Maximum Coverage with Separate Coverage Constraint, with three kinds of solvers: traditional solvers, learning-based solver and specific MILP solver like Gurobi. Code will be made publicly available.

\textbf{4)} The results raise potential concerns about the robustness, such as the $20 \%$ performance decrease for the commercial solver Gurobi. They also help verify the existing theory about data generation~\cite{yehuda2020s} and provide ideas (e.g. adversarial training and parameter tuning) to develop more robust solvers.

\section{Related Works}\vspace{-5pt}
\textbf{Combinatorial optimization solvers.} As a long-standing area, there exist many traditional CO algorithms, including but not limited to greedy algorithms, heuristic algorithms like simulated annealing (SA)~\cite{van1987simulated}, Lin–Kernighan–Helsgaun (LKH3)~\cite{helsgaun2017extension}, and branch-and-bound solvers like CBC~\cite{cbc}, SCIP~\cite{scip}, and Gurobi~\cite{llc2020gurobi}. Besides, driven by the recent development of deep learning and reinforcement learning, many learning-based methods have also been proposed to tackle these problems. A mainstream approach using deep learning is to predict the solution end-to-end, such as the supervised model Pointer Networks~\cite{vinyals2015pointer}, reinforcement learning models S2V-DQN~\cite{khalil2017learning} and MatNet~\cite{kwon2021matrix}. Though these methods did perform well on different types of CO problems, they are not that robust and universal, as discussed in \cite{bengio2020machine}, the solvers may get stuck around poor solutions in many cases. The sensitivity of CO algorithms is theoretically characterized in \cite{varma2021average}, however the metric in \cite{varma2021average} requires optimal solutions which are usually unavailable in practice concerning the NP-hard challenge in CO. \looseness=-1

%\LH{Therefore, the attention should be paid to the hard cases and our framework is proposed to find the cases and utilize them to improve the robustness of the solvers.} \Origin{Different from works~\cite{moon2019parsimonious,Zang_2020} which apply CO for attack against neural networks, we take an initiative on the adversarial attack and defense on CO.}

\textbf{Adversarial attack for neural networks.}
Since the seminal study~\cite{Szegedy2014IntriguingPO} shows that small input perturbations can change model predictions, many adversarial attack methods have been devised to construct such attacks. In general, adversarial attacks can be roughly divided into two categories: white-box attacks with access to the model gradients, e.g.~\cite{goodfellow2014explaining,madry2018towards, Carlini2017TowardsET}, and black-box attacks, with only access to the model predictions, e.g.~\cite{ilyas2018black,narodytska2016simple}. Besides image and text adversarial attacks~\cite{jia2017adversarial}, given the importance of graph-related applications and the successful applications of graph neural networks (GNN)~\cite{scarselli2008graph}, more attention is recently paid to the robustness of GNNs \cite{dai2018adversarial}. 
%In the mean time, many defense strategies like adversarial training~\cite{ganin2016domainadversarial,tramer2020ensemble} have also been proposed to counter this series of attack methods. Since CO problems can usually be encoded by a graph structure and inspired by~\cite{dai2018adversarial}, which develops an RL based attack policy towards GNNs, we propose a novel and flexible attack framework for CO solvers using both heuristic and RL methods.
Besides, the recent adversarial graph matching (GM) network shows how to fulfill attack or defense via perturbing or regularizing geometry property on the GM solver. \cite{gmatt} degrades the quality of GM by perturbing nodes to more dense regions while~\cite{gmdef} improves robustness by separating nodes to be distributed more broadly. However, the techniques are deliberately tailored to the specific problem and can hardly generalize to the general CO problems. \cite{geisler2021generalization} develops a gradient-based attack method for CO learning models that must be differentiable, and cannot generalize to solvers that are usually non-differentiable due to the discrete nature. The black-box differentiation technique~\cite{poganvcic2019differentiation} is restricted for linear objective. We aim to propose a general pipeline by developing more flexible black-box attacks. \looseness=-1

\section{ROCO Framework}\label{sec:robo}\vspace{-5pt}

\subsection{Adversarial Robustness for CO}\label{sec:pre}\vspace{-5pt}

\textbf{CO problem.} In general, a CO problem aims to find the optimal solution under a set of constraints (usually encoded by a latent graph). Typically, we formulate a CO problem $Q$ as:
\begin{equation}
    \min_\mathbf{x}c(\mathbf{x},Q)\ \ s.t.\ \ h_i(\mathbf{x},Q)\leq 0,\ i=1,...,I
\end{equation}
where $\mathbf{x}$ denotes the decision variable (i.e. solution) that should be discrete, $c(\mathbf{x},Q)$ denotes the cost function given problem instance $Q$ and $\{h_i(\mathbf{x},Q)\}_{i=1}^I$ represents the set of constraints. For example, in DAG scheduling, the constraints ensure that the solution $\mathbf{x}$, i.e. the execution order of the DAG job nodes, lies in the feasible space and does not conflict the topological dependency structure of $Q$. However, due to the NP-hard nature, the optimal solution $\mathbf{x}^*$ can be intractable within polynomial time. Therefore, we use a solver $\hat{\mathbf{x}}=f_\theta(Q)$, which gives a mapping $f_\theta: \mathbb{Q} \rightarrow \mathbb{X}$ to \emph{approximate} the optimal solution $\mathbf{x}^*$. In this work, $\theta$ are solvers' parameters (e.g. weights and biases in a neural solver, hyperparameters in Gurobi), $Q \in \mathbb{Q}$ is the problem instance, and $\hat{\mathbf{x}} \in \mathbb{X}$ is the approximated solution. {$\hat{\mathbf{x}}$ also lies in the feasible space that satisfies all the constraints.} 

\textbf{Solvers' robustness.} 
In this paper, we would like to raise the concern that the estimation capability of the solvers can be unstable and fail on certain hard instances. Thus, we need methods to discover the non-trivial hard instances for the solvers as the robustness metric. Here non-trivial means these hard instances are not obtained by trivial operations like increasing the size of problems. In short, given a certain solver $f_\theta$, we design an attacker (i.e.\ perturbation model) $g$ to discover hard instances around existing clean instances $\tilde{Q}=g(f_\theta,Q)$. For example, the attacker can modify the limited number of constraints in $Q$ to get a new problem $\tilde{Q}$. Besides working as a robustness metric, the discovered hard instances can guide the parameter setting of solvers, i.e. design solvers' parameters $\theta$ for better performance (robustness) on hard instances, but it may be beyond the scope of this paper and we mainly focus on developing the first practically feasible robustness metric.

%, both to distinguish more robust solvers and to guide the parameter setting of solvers, i.e. design solvers' parameters $\theta$ for better performance (robustness) on hard instances. Here \emph{practical} ensures that our hard instances should be non-trivial ones, i.e. they are not obtained by trivial operations like increasing the size of problems. In short, given a certain solver $f_\theta$, we design an attacker (i.e. perturbation model) $g$ to discover hard instances around existing clean instances $\tilde{Q}=g(f_\theta,Q)$. For example, the attacker can modify limited number of constraints in $Q$ to get a new problem $\tilde{Q}$.

\textbf{Robustness metric.} A natural evaluation of the solvers' robustness is to use the gap $c(\hat{\mathbf{x}},Q)-c(\mathbf{x}^*,Q)$ \cite{varma2021average}, where a narrower gap stands for a better performance. However, the optimum $\mathbf{x}^*$ is usually unavailable due to the NP-hardness. In this paper, we propose a {robustness evaluation method}, and the attack is defined by:
\begin{definition}[Successful Attack]
Given a certain solver $f_\theta$ and a clean CO problem instance $Q$, we obtain a new problem $\tilde{Q}=g(f_\theta,Q)$ by the attacker $g$ s.t. the optimal cost value after perturbation will become no-worse, i.e. $c(\mathbf{\tilde{x}}^{*},\tilde{Q})\leq c(\mathbf{x}^*,Q)$. A {successful attack} occurs when $c(f_\theta(\tilde{Q}),\tilde{Q}) > c(f_\theta(Q),Q)$.
\end{definition}

\begin{figure}[tb!]
    \centering
    \subfigure[ROCO framework]{
    \label{fig:framework}
    \includegraphics[width=6.5cm]{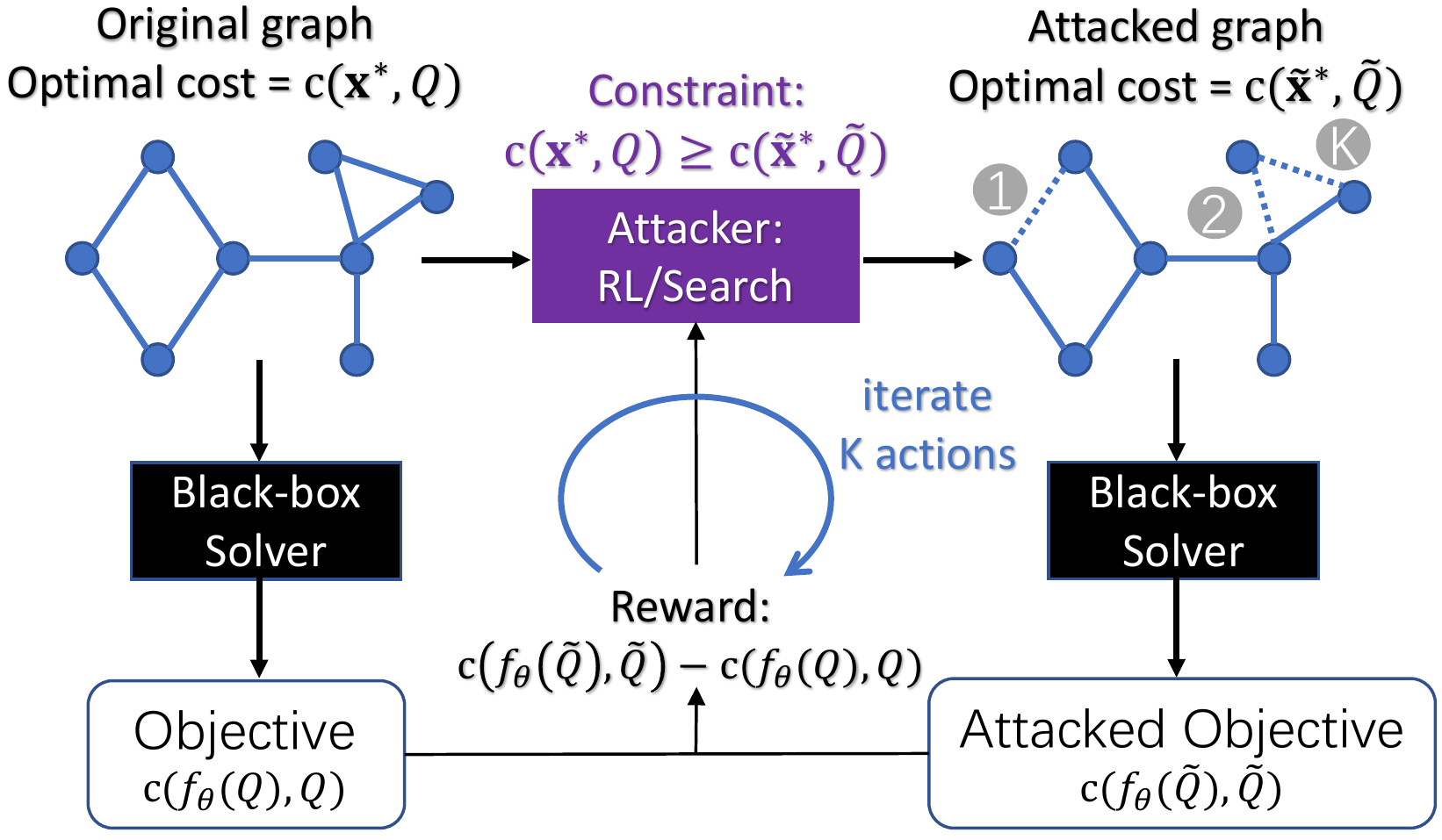}}
    \hspace{0.1in}
    \subfigure[Successful attack]{
    \label{fig:attack}
    \raisebox{0.8cm}{\includegraphics[width=6.5cm]{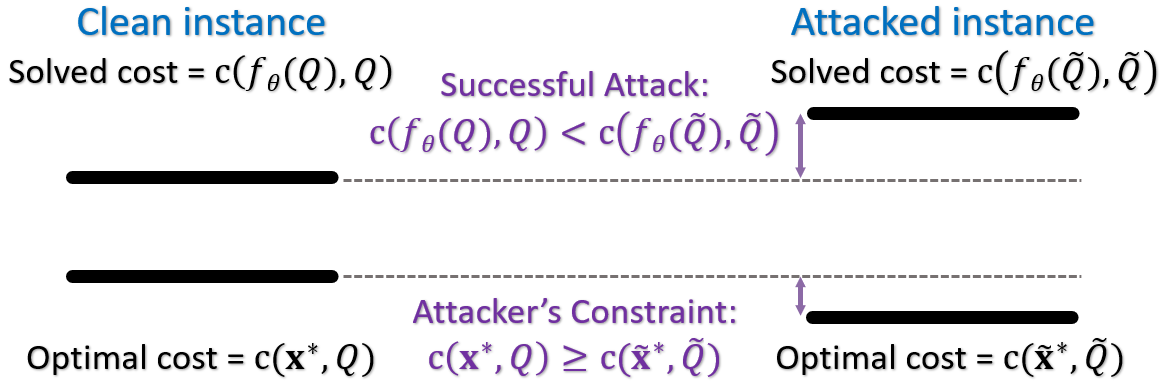}}
    }
    \caption{Overview and intuition of our attack framework. \textbf{Left:} Overview of our attack framework ROCO for CO solvers. ROCO targets on the CO problems which can be encoded by graph (often holds in practice). Here modifying the edges in the encoded graph represents modifying constraints or lowering the cost in CO. \textbf{Right:} The definition of a successful attack.}
    \label{fig:roco}
    \vspace{-10pt}
\end{figure}

The intuition is that we restrict the attacker to discover a new problem instance with no-worse optimal cost value, such that the gap between the solver and the optimal will be definitely enlarged if the solver gives a worse cost value:
\begin{equation}
    \begin{aligned}
    &c(\mathbf{\tilde{x}}^{*},\tilde{Q})\leq c(\mathbf{x}^*,Q),\ c(f_\theta(\tilde{Q}),\tilde{Q}) > c(f_\theta(Q),Q) \Rightarrow  \\
    &c(f_\theta(Q),Q)-c(\mathbf{x}^*,Q) < c(f_\theta(\tilde{Q}),\tilde{Q})-c(\mathbf{\tilde{x}}^*,\tilde{Q}),
    \end{aligned}
\end{equation}
which represents a successful attack (see Fig.~\ref{fig:attack}).

\textbf{Attacker's target.} In conclusion, the attacker should maximize the solver's cost value after perturbation while conforming to the \emph{no-worse} restriction. Therefore, we formulate our attack objective:
\begin{equation}\label{eq:attack}
\begin{aligned}
    &\max_{\tilde{Q}}\quad c(f_\theta(\tilde{Q}),\tilde{Q})-c(f_\theta(Q),Q) \\
    &s.t.\quad \tilde{Q}=g(f_\theta,Q),\ c(\mathbf{\tilde{x}}^{*},\tilde{Q})\leq c(\mathbf{x}^*,Q)
\end{aligned}
\end{equation}

\textbf{Attacker's constraints.} Two types of strategies are developed to meet the \emph{no worse optimal cost} constraint for the attackers: \textbf{1)}~Loosening the constraints of the problem. In this way, the original optimal solution still lies in the feasible space and its cost is not affected, while the enlarged feasible space may include better solutions. Therefore, the new optimal cost cannot be worse than the original one, and we discover this strategy effective for DAG scheduling (see Section~\ref{sec:DAG}). \textbf{2)}~Modifying a partial problem by lowering the costs. Such a strategy ensures that the cost that corresponds to the original optimal solution can never become worse, and there are also chances to introduce solutions with lower costs. We demonstrate its successful application in traveling salesman (see Section~\ref{sec:ATSP}) and max coverage problems (see Section~\ref{sec:MC} \& Appendix \ref{app:experimentMCSCC}).

\textbf{Attack via graph modification.}
As we have mentioned above, a CO problem $Q$ can usually be encoded by a latent graph $\mathcal{G}=(V, E)$. For example, the DAG problem can be represented by a directed acyclic graph (see Fig.~\ref{fig:DAGexample}), s.t. the graph nodes stand for the DAG jobs, the directed edges denote the sequential constraints for the jobs, and the node features represent the jobs' completion time. Hence, the perturbation to the problem instance can be conducted by modifications to the graph structure~\citep{dai2018adversarial}. Besides, to get non-trivial hard instances around clean examples, the attack is limited within a given attack budget $K$, e.g. the attacker can remove $\leq K$ edges in a DAG problem $Q$. Typically, in our setting, the graph modification is all conducted by \emph{modifying edges} in the graph.

%In Section~\ref{sec:MDP}, we devise both reinforce learning (RL) and heuristic based attackers to discover hard instances from the clean test examples. For RL, the popular Proximal Policy Optimization (PPO)~\cite{schulman2017proximal} framework is adopted. Besides, we also implement three other traditional attack algorithms: random sampling, optimum-guided search and simulated annealing.
\subsection{Implementations of the Attacker}\label{sec:MDP}\vspace{-5pt}
\textbf{Random attack baseline.} The baseline simply chooses and modifies $K$ random edges sequentially, then calculates the new solution cost. The baseline is designed to reflect the instability of CO problems, we compare other attackers (which achieve much higher attack effects) to the random attack baseline to claim that \emph{our framework can uncover far more hard instances than the inherent instability caused by the CO problems}.  

\textbf{Reinforcement learning attacker} \textit{(ROCO-RL)}.
We use Eq.~\ref{eq:attack} as the objective and resort to reinforcement learning (RL) to find $\tilde{Q}$ in a data-driven manner. Specifically, we modify the graph structure and compute $c(f_\theta(\tilde{Q}),\tilde{Q})$ alternatively, getting rewards that will be fed into the PPO~\cite{schulman2017proximal} framework and train the agent.

Given an instance ($f_\theta,Q$) along with a modification budget $K$, we model the attack via sequential edge modification as a Markov Decision Process (MDP):

\textit{State.} The current problem instance $Q^k$ (i.e. the problem instance $\tilde{Q}$ obtained after $k$ actions) is treated as the state. The original problem instance $Q^0=Q$ is the starting state.

\textit{Action.} The attacker is allowed to modify edges in the graph. So a single action at time step $k$ is $a^k \in \mathcal{A}^k \subseteq E^k$. Here our action space $\mathcal{A}^k$ is usually a subset of all the edges $E^k$ since we can narrow the action space (i.e. abandon some useless edge candidates) according to the previous solution $f_\theta(Q^{k})$ to speed up our algorithm. Furthermore, we decompose the action space ($O(|V|^2) \rightarrow O(|V|)$) by transforming the edge selection into two sequential node selections: first selecting the start node, then the end node. \looseness=-1

\textit{Reward.}  The objective of the new CO problem $Q^{k + 1}$ is $c(f_\theta(Q^{k+1}),Q^{k+1})$. To maximize the new cost value, we define the reward as the increase of the objective:
\begin{equation}
    r = c(f_\theta(Q^{k+1}),Q^{k+1}) - c(f_\theta(Q^{k}),Q^{k})
\end{equation}

\textit{Terminal.} Once the agent modifies $K$ edges or edge candidates become empty, the process stops. 

The input and constraints of a CO problem can often be encoded in a graph $\mathcal{G}=(V,E)$, and the PPO agent (i.e. the actor and the critic) should behave according to the graph features. We resort to the Graph Neural Networks (GNNs) for graph embedding:
\begin{equation}\label{GraphEmbed}
    \textbf{n} = \mbox{GNN}(\mathcal{G}^k),\ \textbf{g}=\mbox{AttPool}(\textbf{n})
\end{equation}
where $\mathcal{G}^k$ is the corresponding graph of problem $Q^k$, the matrix $\textbf{n}$ (with the size of node number $\times$ embedding dim) is the node embedding. An attention pooling is used to extract a graph level embedding $\textbf{g}$. The GNN model can differ by the CO problems (details in Appendix~\ref{app:graph_embedding}). After graph feature extraction, we design the actor and critic nets:

\textit{Critic.} The critic predicts the value of each state $Q^k$. Since it aims at reward maximization, a max-pooling layer is adopted over all node features which are concatenated (denoted by $[\cdot || \cdot]$) with the graph embedding $\mathbf{g}$, fed into a network (e.g. ResNet block~\cite{he2016deep}) for value prediction:
\begin{equation}\label{Value}
    \mathcal{V}(Q^k) = \mbox{ResNet}_1([\mbox{MaxPool}(\textbf{n})||\textbf{g}])
\end{equation}

\textit{Actor.} The edge selection is implemented by selecting the start and end node sequentially. The action scores are computed using two independent ResNet blocks, and a Softmax layer is added to regularize the scores into probability $[0,1]$:
\begin{equation}\label{Prob}
\begin{aligned}
    P(a_1) &= \mbox{Softmax}(\mbox{ResNet}_2([\textbf{n}||\textbf{g}])),\\ 
    P(a_2|a_1) &= \mbox{Softmax}(\mbox{ResNet}_3([\textbf{n}||\textbf{n}[a_1]||\textbf{g}))
\end{aligned}
\end{equation}
where $\textbf{n}[a_1]$ denotes node $a_1$'s embedding. We add the feature vector of the selected start node for the end node selection. For training, actions are sampled w.r.t. their probabilities. For testing, beam search is used to find the optimal solution: actions with top-$B$ probabilities are chosen for each graph in the last time step, and only those actions with top-$B$ rewards will be reserved for the next search step (see Alg.~\ref{alg:PPO}). \looseness=-1

\begin{wrapfigure}{r}{0.55\textwidth}
 \begin{minipage}[t]{1\linewidth}
 \centering
 \vspace{-3pt}
\begin{algorithm2e}[H]
{\small{
    \caption{\textbf{Attack by iterative edge manipulation}\label{alg:PPO}}
    \KwIn{Input problem $Q$; solver $f_\theta$; max number of actions $K$; beam size $B$.}% bias $eps$
    $Q^0_{1..B} \leftarrow Q$; $\tilde{Q} \leftarrow Q$; \algcomment{set initial state}\\
    \For {$k \leftarrow 1..K$} 
    {
        \For {$b \leftarrow 1..B$}{ \algcomment{do beam search for problems in last step} \\
            Predict $P(a_1)$, $P(a_2|a_1)$ on $Q^{k-1}_b$;\\
            Select $(a_1,a_2)$ with top-$B$ probabilities;  
        }
        \For{\rm{each} $(b,a_1,a_2)$ pair}{
            $Q^\prime(b,a_1,a_2) \leftarrow $ modify edge $(a_1,a_2)$ in $Q^{k-1}_b$;  \\
             \If{$c(f_\theta,Q^\prime(b,a_1,a_2)) > c(f_\theta,Q^*)$}{
                $\tilde{Q} \leftarrow Q^\prime(b,a_1,a_2)$ \\ \algcomment{update the best attacked problem}
            }
        }
        Sort $Q^\prime(\cdot,\cdot,\cdot)$ w.r.t. cost values by decreasing order;\\ 
        $Q^k_{1..B} \leftarrow Q^\prime_{1..B}$; \algcomment{select top-$B$ problems for the next step}
    }
    \KwOut{Best attacked problem instance $\tilde{Q}$.}}
    }
\end{algorithm2e}
 \end{minipage}\vspace{-20pt}
 \end{wrapfigure}

\textbf{Other attacker implementations.}
We also implement three other traditional attack algorithms for comparison: random search, optimum-guided search and simulated annealing.

\emph{i) Random Search (ROCO-RA):}
In each iteration, an edge is randomly chosen to be modified in the graph and it will be repeated for $K$ iterations.  Different from the random attack baseline without any search process, we will run $N$ attack trials and choose the best solution. Despite its simplicity, it can reflect the robustness of solvers with the cost of time complexity $O(NK)$. 

\emph{ii) Optimum-Guided Search (ROCO-OG):}
This method focuses on finding the optimum solution during each iteration. We use beam search to maintain the best $B$ current states and randomly sample $M$ different actions from the candidates to generate next states. The number of iterations is set to be no more than $K$, so its time complexity is $O(BMK)$.

\emph{iii) Simulated Annealing (ROCO-SA):} Simulated annealing \cite{van1987simulated} comes from the idea of annealing and cooling used in physics for particle crystallization. In our scenario, a higher temperature indicates a higher probability of accepting a worse solution, allowing to jump out of the local optimum. As the action number increases, the temperature will decrease, and we are tending to reject the bad solution. The detailed process is shown in Appendix~\ref{app:SA} and we will repeat the algorithm for $N$ times. SA is a fine-tuned algorithm and we can use grid search to find the best parameter to fit the training set. Its time complexity is $O(NMK)$. 

Tab.~\ref{tab:timecomp} in Appendix \ref{app:exp} concludes the attack methods property and time complexity. Since these three traditional algorithms are inherently stochastic, we run them multiple times to calculate the mean and standard deviation. Meanwhile, the reinforcement learning attacker is deterministic because we always use the best trained model to modify graph structures on the testing set.

% \subsection{Hard Case Utilization and Defense}
% Apart from the robustness evaluation of combinatorial optimization solvers, the attackers can help generate the hard cases for the given data samples. A better utilization of the cases can help improve the robustness or even performance of the solvers under the similar data distribution. The process requires the solvers can be tuned rather than fixed. Therefore, we take the neural network model and the Gurobi solver as the examples to verify the usability of the hard cases. Adversarial training and parameter tuning are applied and more details can be seen in section~\ref{sec:exp}.

\section{Studied CO Tasks and Results} \label{sec:exp}\vspace{-5pt}
We conduct experiments on four representative CO tasks: Directed Acyclic Graph Scheduling, Asymmetric Traveling Salesman Problem, Maximum Coverage and Maximum Coverage with Separate Coverage Constraint. Considering the page limitation, the details of the last task is postponed to Appendix~\ref{app:experimentMCSCC}.
We also test the robustness of representative solvers including heuristics~\cite{helsgaun2017extension,grandl2014multi}, branch-and-bound solvers~\cite{cbc,scip,llc2020gurobi} and learning-based methods~\cite{kwon2021matrix}. The detailed graph embedding methods for the four tasks is shown in Appendix~\ref{app:graph_embedding}. In Appendix \ref{app:exp}, we provide training and evaluation hyperparameters of different solvers for fair time comparison and reproducibility. Experiments are run on our heterogeneous cluster with RTX 2080Ti and RTX 3090.

\subsection{Task I: Directed Acyclic Graph Scheduling} \label{sec:DAG}
\begin{wrapfigure}{r}{0.4\textwidth}
\begin{minipage}[t]{1\linewidth}
% \begin{figure}[h]
    \centering
    \includegraphics[width=\linewidth]{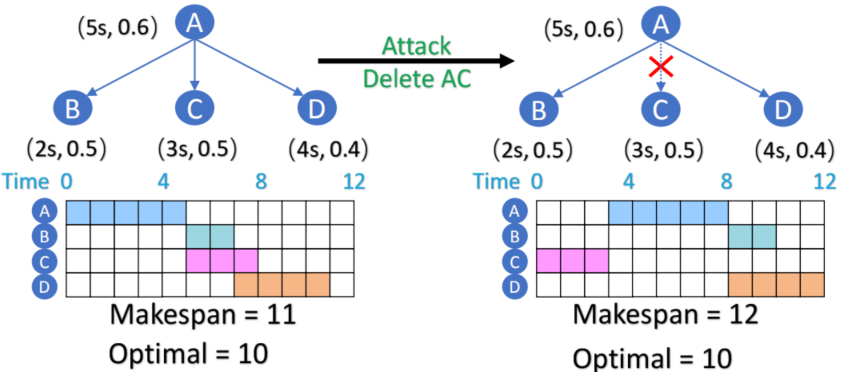} \vspace{-15pt}
    \caption{DAG attack on Shortest Job First algorithm. The edges show  sequential dependencies. The $(x, y)$ tuple of each node means running time ($x$) and resource occupancy rate ($y$).}% Similar sketches for ATSP and MCSCC (see experiments) are in Appendix~\ref{app:examples} 
    \label{fig:DAGexample}   \vspace{-19pt}
% \end{figure}
\end{minipage}  
\end{wrapfigure}

\vspace{-5pt}
Task scheduling for heterogeneous systems and various jobs is a practical problem. Many systems formulate the job stages and their dependencies as a Directed Acyclic Graph (DAG)~\citep{saha2015apache, chambers2010flumejava, zaharia2012resilient}, as shown in Fig.~\ref{fig:DAGexample}. The data center has limited computing resources to allocate the jobs with different resource requirements. These jobs can run in parallel if all their parent jobs have finished and the required resources are available.  Our goal is to \textbf{minimize the makespan (i.e.\ finish all jobs ASAP.)}

\begin{table*}[tb!]
    \centering%TPC-H-$X$ means $X$ jobs are jointly scheduled. 
    \caption{Evaluation of solvers for DAG scheduling. Clean makespans are reported on the clean test set.  Larger attacked makespans mean stronger attacker. Random denotes the random attack baseline with the mean and std tested for 100 trials. ROCO-RA/OG/SA are tested for 10 trials to calculate the mean and std. ROCO-RL is reported with a single value since it is invariant to different random seeds. All four attack methods perform better than the random attack baseline, which usually leads to a shorter makespan (unsuccessful attack), proving that our attack framework can find more hard instances than those because of the inherent instability of CO problems. Tetris is the most robust solver on dataset sizes of 50 and 150 while Critical path is more robust on size 100.}
    \resizebox{\textwidth}{!}
    {
    \begin{tabular}{r|c|c|ccccc}
        \hline
        % Table generated by Excel2LaTeX from sheet 'MGM-willow'
        \multirow{2}[0]{*}{Solver} & {Problem Size:} & Clean Makespan &  \multicolumn{5}{c}{Attack Method (Attacked Makespans  $(\times 10^5)$)}  \\
         & \#job & $(\times 10^5) \downarrow$ & Random & ROCO-RA & ROCO-OG & ROCO-SA & ROCO-RL   \\
        \hline
        Shortest Job First & 50  & 1.0461 & $1.0406 \pm 0.0035$  & $1.0574 \pm 0.0013$  & $1.0600 \pm 0.0019$ & $\textbf{1.0622} \pm \textbf{0.0007}$ & $1.0608$   \\
        Critical Path & 50  & 0.8695 & $0.8828 \pm 0.0077$ & $0.9402 \pm 0.0038$ & $0.9480 \pm 0.0021$  & $\textbf{0.9528} \pm \textbf{0.0013}$ &  $0.9500$  \\
        Tetris~\citep{grandl2014multi} & 50  & 0.8227  & $0.8582 \pm 0.0079$ & $0.0952 \pm 0.0049$ & $0.9218\pm 0.0066$ & $0.9380 \pm 0.1153$ & $\textbf{0.9397}$   \\
        \hline
        Shortest Job First & 100  & 1.9160 & $1.8832 \pm 0.0075$ & $1.9210 \pm 0.0006$ & $1.9239 \pm 0.0008 $ & $1.9239 \pm 0.0008$   &  $\textbf{1.9263}$  \\
        Critical Path & 100  & 1.6018 & $1.6279 \pm 0.0167$ & $1.7391 \pm 0.0045$ & $1.7456 \pm 0.0043$ & $1.7480 \pm 0.0003$ &  $\textbf{1.7498}$ \\
        Tetris~\citep{grandl2014multi} & 100  & 1.5186  & $1.5830 \pm 0.0196$ & $1.7201 \pm 0.0055$ & $1.7099 \pm 0.0111$ & $1.7418 \pm 0.0074 $  & $\textbf{1.7526}$  \\
        \hline
        Shortest Job First & 150 & 2.8578 & $2.8281 \pm 0.0074$ & $2.8818 \pm 0.0020$ & $2.8898 \pm 0.0023$  & $2.8950 \pm 0.0014$  &  $\textbf{2.8964}$  \\
        Critical Path & 150  & 2.4398 & $2.4254 \pm 0.0203$ & $2.5698 \pm 0.0093$ & $2.5928 \pm 0.0090$ &  $2.6020 \pm 0.0029$  & $\textbf{2.6069}$  \\
        Tetris~\citep{grandl2014multi} & 150  & 2.2469  & $2.3581 \pm 0.0270$ & $2.4988 \pm 0.0191$ & $2.5039 \pm 0.0202$ & $\textbf{2.5399} \pm \textbf{0.0058}$ &  $2.5329$ \\
        \hline
    \end{tabular}%
    }
    \vspace{-22pt}
    \label{tab:DAG_attack}
\end{table*}

\textbf{Solvers.} We choose three popular heuristic solvers as our attack targets. First, the {Shortest Job First} algorithm chooses the jobs greedily with minimum completion time. Second, the~{Critical Path} algorithm finds the bottlenecks and finishes the jobs in the critical path sequence. Third, the {Tetris} \citep{grandl2014multi} scheduling algorithm models the jobs as 2-dimension blocks in the Tetris games according to their finishing time and resource requirement. 

\textbf{Attack model.} The edges in a DAG represent job dependencies, and removing edges will relax the constraints. After removing existing edges in a DAG, it is obvious that the new solution will be equal to or better than the original one since there are fewer restrictions. As a result, in the DAG scheduling tasks, the attack model is to selectively remove existing edges.

\textbf{Dataset.} We use TPC-H dataset~(\url{http://tpc.org/tpch/default5.asp}), composed of business-oriented queries and concurrent data modification. Many DAGs referring to computation jobs, have tens or hundreds of stages with different duration and numbers of parallel tasks. We gather the DAGs randomly and generate three different sizes of datasets, TPC-H-50, TPC-H-100, TPC-H-150, each with $50$ training and $10$ testing samples. The DAG nodes have two properties: execution time and resource requirement. 

\begin{wrapfigure}{r}{0.50\textwidth}
\begin{minipage}[t]{1\linewidth}
    \centering
    \includegraphics[width=0.48\linewidth]{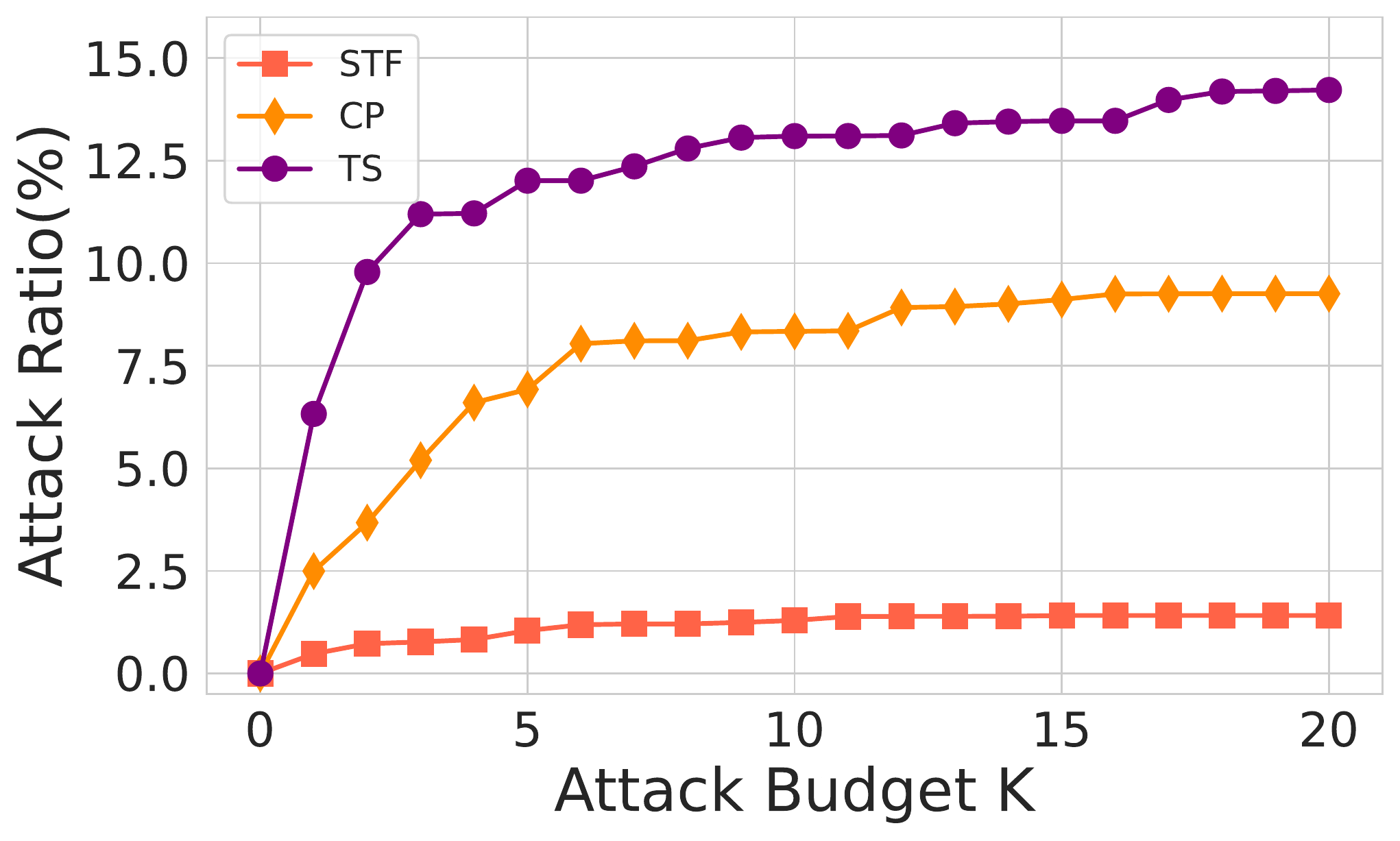}
    \includegraphics[width=0.48\linewidth]{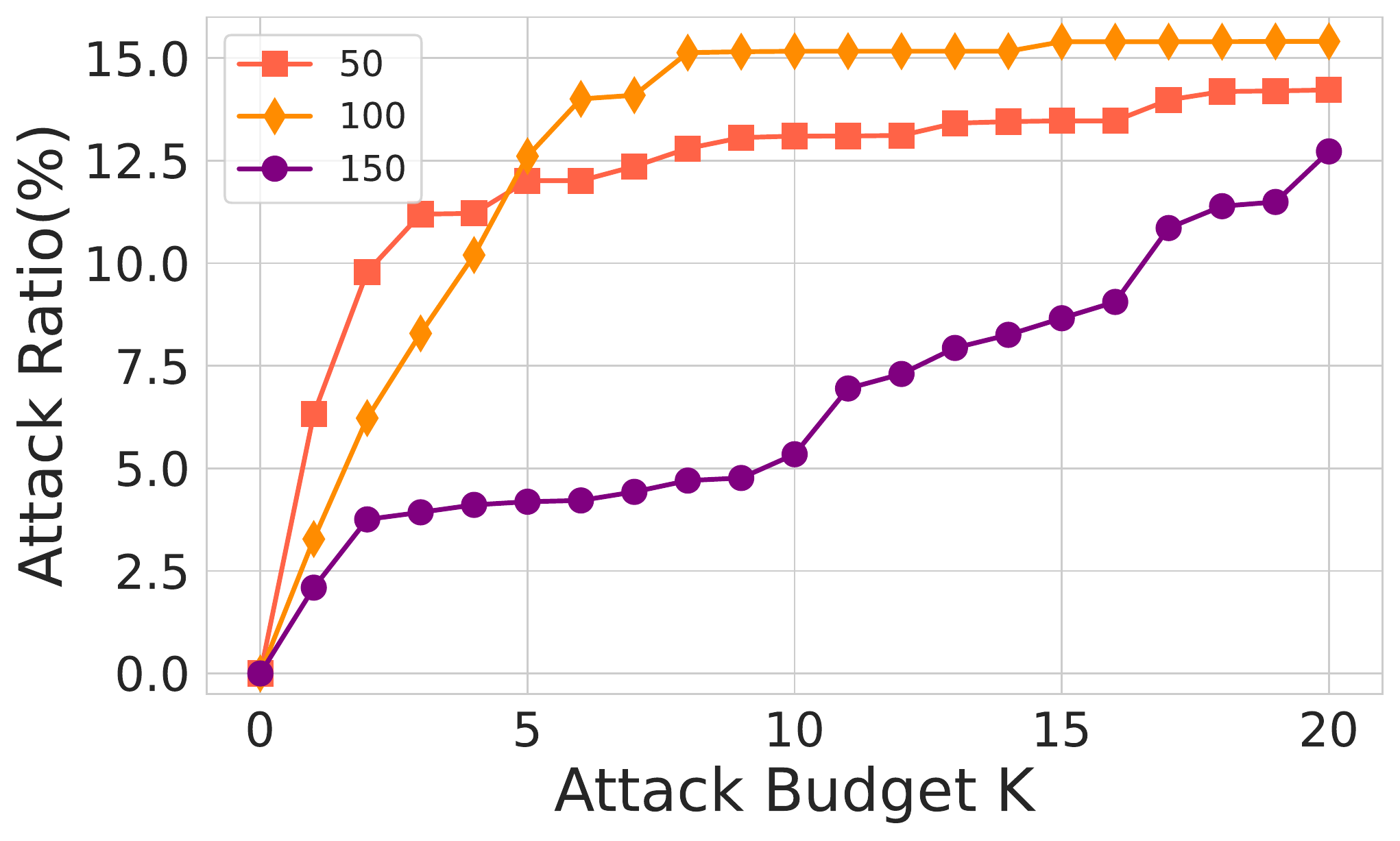}
    \vspace{-10pt}
    \caption{Left: RL Attack results w.r.t. attack budgets $K$ on TPC-H-50 dataset for solvers Shortest Job First(SJF), Critical Path(CP) and Tetris(TS). Ratio represents the increase of makespan after attack w.r.t.\ clean makespans. Right: RL Attack results w.r.t.\ attack budgets $K$ on TPC-H-50/100/150 datasets for the {Tetris (TS)} solver.}
    \label{fig:DAGEdges} 
\end{minipage}  
\vspace{-10pt}
\end{wrapfigure}

\textbf{Results and analysis.} Tab.~\ref{tab:DAG_attack} reports the results of our four attack methods.
In the designed attack methods, RL outperforms other learning-free methods in most cases, illustrating the correctness of our feature extraction techniques and training framework. It is worth noting that even the simplest random search (ROCO-RA) can cause a significant performance degradation to the CO solvers, showing their vulnerability and the effectiveness of the attack framework. Fig.~\ref{fig:DAGEdges} (left) shows the attack ratio for different solvers on TPC-H-50. In general, the attack ratio increases then flattens out w.r.t. attack budgets, corresponding with the fact that we can find more hard instances under a larger searching space (but harder to find new successfully attacked instances as we remove more constraints). Fig.~\ref{fig:DAGEdges} (right) demonstrates the effect of attack budgets on different sizes of datasets using the same solver Tetris. The figure illustrates that attack ratios on smaller datasets tend to flatten first (larger datasets allow for more attack actions), which can guide us to select the right number of attack actions for different datasets. \looseness=-1

\subsection{Task II: Asymmetric Traveling Salesman } \label{sec:ATSP}\vspace{-5pt}
The classic traveling salesman problem (TSP) is to \textbf{find the shortest cycle to travel across all the cities.} Here we tackle the extended version namely asymmetric TSP (ATSP) for its generality.

%\subsubsection{Implementation Details}\vspace{-8pt}
\textbf{Solvers.} The robustness of four algorithms are evaluated: First, {Nearest Neighbour} greedily adds the nearest city to the tour. Second, {Furthest Insertion} finds the city with the furthest distance to the existing cities in the tour and inserts it. Third, {Lin-Kernighan Heuristic (LKH3)}~\citep{helsgaun2017extension} is the traditional SOTA TSP solver. Finally, {Matrix Encoding Networks (MatNet)}~\citep{kwon2021matrix} claims as a SOTA learning-based solver for ATSP and flexible flow shop (FFSP).

\textbf{Attack model.} The attack is to choose an edge and half its value, after which we will get a no-worse theoretical optimum since the length of any path will maintain the same or even decrease. To reduce the action space, we will not select the edges in the current path predicted by the solver at the last time step. A successful attack example is shown in Fig.~\ref{fig:TSPexample} in Appendix \ref{app:examples}.

\begin{table*}[tb!]
    \centering
    \caption{Evaluation of solvers for ATSP. 
    A Larger tour length means a stronger attacker. Random denotes the random attack baseline with the mean and std tested for 100 trials. ROCO-RA/OG/SA are tested for 10 trials to calculate the mean and std. ROCO-RL is reported with a single value since it is invariant to different random seeds. MatNet (fixed) represents training MatNet with 10000 fixed ATSP instances. LKH3 and MatNet are more robust than peer solvers.}
    \resizebox{\textwidth}{!}
    {
    \begin{tabular}{r|c|c|ccccc}
        \hline
        % Table generated by Excel2LaTeX from sheet 'MGM-willow'
        \multirow{2}[0]{*}{Solver} & Problem  & Clean Tour   & \multicolumn{5}{c}{Attack Method (Attacked Tour Length $Length (\times 10^6)$)}  \\
         &Size: \# city & $(\times 10^6) \downarrow$ & Random & ROCO-RA & ROCO-OG & ROCO-SA & ROCO-RL   \\
        \hline
       
        Nearest Neighbour & 20  & 1.9354 & $ 1.9143 \pm 0.0259 $ & $ 2.1307 \pm 0.0153 $ & $ 2.1162 \pm 0.0323 $ & $ 2.1344 \pm 0.0159 $ & $ \textbf{2.1858}$   \\
        Furthest Insertion & 20  & 1.6092  & $ 1.5847 \pm 0.0142 $ & $ 1.6953 \pm 0.0105 $ & $ 1.6926 \pm 0.0117 $ & $ 1.7183 \pm 0.0114 $ &  $ \textbf{1.7469}$  \\
        LKH3~\citep{helsgaun2017extension} & 20  & 1.4595 & $ 1.4201 \pm 0.0089 $ & $ 1.4599 \pm 0.0003 $ & $ 1.4599 \pm 0.0004 $ & $ 1.4610 \pm 0.0017 $ &  $\textbf{1.4611}$ \\
        MatNet~\citep{kwon2021matrix} & 20 & 1.4616 & $ 1.4237 \pm 0.0089 $ & $ 1.4674 \pm 0.0012 $ & $ 1.4683 \pm 0.0006 $ & $ 1.4683 \pm 0.0009 $ &  $\textbf{1.4708}$ \\
        MatNet (fixed)~\citep{kwon2021matrix} & 20 & 1.4623 & $ 1.4227 \pm 0.0086 $ & $ 1.4665 \pm 0.0013 $ & $ 1.4687 \pm 0.0004 $ & $ 1.4714 \pm 0.0013 $ & \textbf{1.4737} \\
        \hline
        
        Nearest Neighbour & 50  & 2.2247 & $ 2.2221 \pm 0.0131 $ & $ 2.3635 \pm 0.0082 $ & $ 2.3809 \pm 0.0096 $ & $ 2.4058 \pm 0.0151 $ & $ \textbf{2.4530}$   \\
        Furthest Insertion  & 50  & 1.9772 & $ 1.9730 \pm 0.0091 $ & $ 2.0593 \pm 0.0071 $ & $ 2.0466 \pm 0.0125 $ & $ 2.0632 \pm 0.0089 $ &  $ \textbf{2.1150}$  \\
        LKH3~\citep{helsgaun2017extension} & 50  & 1.6621  & $ 1.6493 \pm 0.0042 $ & $ 1.6653 \pm 0.0007 $ & $ 1.6656 \pm 0.0007 $ & $ \textbf{1.6682} \pm \textbf{0.0010} $  &  $1.6679$ \\
        MatNet~\citep{kwon2021matrix} & 50 & 1.6915 & $ 1.6786 \pm 0.0051 $ & $ 1.7153 \pm 0.0012 $ & $ 1.7204 \pm 0.0016 $ & $ 1.7255 \pm 0.0012 $ &  $\textbf{1.7279}$ \\
        MatNet(fixed)~\citep{kwon2021matrix} & 50 & 1.6972 & $ 1.6894 \pm 0.0054 $ & $ 1.7171 \pm 0.0021 $ & $ 1.7205 \pm 0.0005 $ & $ 1.7279 \pm 0.0007 $ & \textbf{1.7349} \\
        
        \hline
        Nearest Neighbour & 100  & 2.1456 & $ 2.1423 \pm 0.0079 $ & $ 2.2319 \pm 0.0099 $ & $ 2.2213 \pm 0.0152 $ & $ 2.2273 \pm 0.0112 $ & $ \textbf{2.2533}$   \\
        Furthest Insertion  & 100  & 1.9209 & $ 1.9188 \pm 0.0061 $ & $ 1.9762 \pm 0.0088 $ & $ 1.9782 \pm 0.0111 $ & $ 1.9853 \pm 0.0063 $ &  $ \textbf{2.0144}$  \\
        LKH3~\citep{helsgaun2017extension} & 100  & 1.5763 & $ 1.5716 \pm 0.0016 $ & $ 1.5826 \pm 0.0006 $ & $ 1.5848 \pm 0.0005 $ & $ 1.5856 \pm 0.0003 $  &  $\textbf{1.5862}$ \\
        MatNet~\citep{kwon2021matrix} & 100 & 1.6545 & $ 1.6492 \pm 0.0036 $ & $ 1.6772 \pm 0.0012 $ & $ 1.6815 \pm 0.0005 $ & $ 1.6841 \pm 0.0007 $ &  $\textbf{1.6873}$ \\
        MatNet (fixed)~\citep{kwon2021matrix} & 100 & 1.6556 & $ 1.6498 \pm 0.0037 $ & $ 1.6788 \pm 0.0008 $ & $ 1.6823 \pm 0.0003 $ & $ 1.6846 \pm 0.0015 $ & \textbf{1.6919}\\
        \hline
    \end{tabular}%
    }
    \vspace{-15pt}
    \label{tab:ATSP_attack}
\end{table*}

% \textbf{4) Graph embedding.}
%\subsubsection{Experiment results} \vspace{-8pt}
\textbf{Dataset.} The dataset comes from \citep{kwon2021matrix}, consisting of `tmat' class ATSP instances which have the triangle inequality and are widely studied by the operation research community~\citep{cirasella2001asymmetric}. We solve the ATSP of three sizes, $20$, $50$ and $100$ cities.
% The dataset is TSPLIB95 format~\citep{reinelt1995tsplib95} and its 
The distance matrix is fully connected and asymmetric, and each dataset consists of $50$ training samples and $20$ testing samples. 

\textbf{Results and analysis.} Tab.~\ref{tab:ATSP_attack} reports the attack results of four target solvers.  
where the {RL} based attack outperforms other methods in most cases. We can also conclude that LKH3 is the most robust among all compared solvers, possibly because its local search nature can jump out of local optimums. MatNet generates different random problem instances for training in each iteration. The comparison between ``MatNet'' and ``MatNet (fixed)'' proves that the huge amount of i.i.d. training instances promotes both the performance and the robustness of the model.
As~\cite{yehuda2020s} points out, learning-based solvers for CO cannot be trained sufficiently due to the theoretical limitation of data generation. Thus, the training paradigm with large unlabeled data can be a good choice for learning-based CO solvers.

\subsection{Task III: Maximum Coverage} \label{sec:MC}\vspace{-5pt}
The maximum coverage (MC) problem is a classical problem that is widely studied in approximation algorithms. Specifically, as input we will be given several sets $\mathcal{S}=\{S_1,S_2,...,S_m\}$ and a number $k$. The sets consist of weighted elements $e \in \mathcal{E}$ and may have overlap with each other. Our goal is to \textbf{select at most $k$ of these sets such that the maximum weight of elements are covered.} Formally, the above problem can be formulated as follows:
\begin{equation}
\label{eq:mc}
    \max_{\mathcal{S}'}\sum_{e\in\mathcal{E}}w(e)\times\mathbb{I}(e\in\bigcup_{S_i \in \mathcal{S}'}S_i) \quad s.t. \quad \mathcal{S}'\subseteq\mathcal{S},\ |\mathcal{S}'|\leq k
\end{equation}
where $w(\cdot):\mathcal{E}\rightarrow\mathbb{R}$ denotes the weight of a certain element.

\begin{table*}[tb!]
    \centering
    \caption{Evaluation of solvers for MC. The time limit is shown in brackets for branch-and-bound solvers, which should be long enough to output a feasible solution in the clean dataset.   Smaller attacked weight means stronger attacker. Random denotes the random attack baseline with the mean and std tested for 100 trials. ROCO-RA/OG/SA are tested for 10 trials to calculate the mean and std. ROCO-RL is reported with a single value since it is invariant to different random seeds. Gurobi turns out to be the most robust model among all compared solvers, possibly because of its SOTA branch-and-bound strategy.}\label{tab:MCAttack}
    \resizebox{\textwidth}{!}
    {
    \begin{tabular}{r|c|c|ccccc}
        \hline
        % Table generated by Excel2LaTeX from sheet 'MGM-willow'
        \multirow{2}[0]{*}{Solver} & Problem Size: & Clean Weight & \multicolumn{5}{c}{Attack Method (Attacked Weight $(\times 10^4)$)}  \\
         &  \#set-\#element & $(\times 10^4) \uparrow$ & Random & ROCO-RA & ROCO-OG & ROCO-SA & ROCO-RL   \\
        \hline
        Greedy & 100-200 & 0.7349 & $ 0.7353 \pm 0.0007 $ & $ 0.7252 \pm 0.0012 $ & $ 0.7255 \pm 0.0013 $ & $ 0.7253 \pm 0.0007 $ & \textbf{0.7239} \\
        CBC(0.5s)  & 100-200  & 0.7184 & $ 0.7020 \pm 0.0122 $ & $ 0.6161 \pm 0.0206 $ & $ 0.5448 \pm 0.0255 $ & $ 0.5423 \pm 0.0315 $ & \textbf{0.5029}  \\
        SCIP(0.5s) & 100-200  & 0.7115 & $ 0.7116 \pm 0.0053 $ & $ 0.6850 \pm 0.0013 $ & $ 0.6808 \pm 0.0025 $ & $ 0.6803 \pm 0.0015 $ & \textbf{0.6788}  \\
        Gurobi(0.5s) & 100-200 & 0.7491 & $ 0.7495 \pm 0.0004 $ & $ 0.7460 \pm 0.0004 $ & $ 0.7441 \pm 0.0006 $ & $ 0.7426 \pm 0.0008 $ & \textbf{0.7409} \\
        \hline
        Greedy & 150-300 & 1.1328 & $ 1.1330 \pm 0.0011 $ & $ 1.1207 \pm 0.0019 $ & $ 1.1183 \pm 0.0021 $ & $ 1.1201 \pm 0.0017 $ & \textbf{1.1182}   \\
        CBC(1s)  & 150-300 & 1.0510 & $ 1.0181 \pm 0.0192 $ & $ 0.8975 \pm 0.0386 $ & $ 0.8054 \pm 0.0487 $ & $ 0.7959 \pm 0.0482 $ & \textbf{0.7527}  \\
        SCIP(1s) & 150-300 & 1.0847 & $ 1.0867 \pm 0.0027 $ & $ 1.0544 \pm 0.0017 $ & $ 1.0467 \pm 0.0016 $ & $ 1.0486 \pm 0.0012 $ & \textbf{1.0443} \\
        Gurobi(1s) & 150-300 & 1.1485 & $ 1.1486 \pm 0.0009 $ & $ 1.1391 \pm 0.0008 $ & $ 1.1358 \pm 0.0011 $ & $ 1.1332 \pm 0.0005 $ & \textbf{1.1314}\\
        \hline
        Greedy & 200-400 & 1.4659 & $ 1.4665 \pm 0.0013 $ & $ 1.4537 \pm 0.0019 $ & $ 1.4533 \pm 0.0019 $ & $ 1.4552 \pm 0.0013 $ & \textbf{1.4466}   \\
        CBC(1.5s)  & 200-400 & 1.4248 & $ 1.4088 \pm 0.0219 $ & $ 1.2053 \pm 0.0475 $ & $ 1.0966 \pm 0.0575 $ & $ 1.0866 \pm 0.0705 $ & \textbf{1.0418}  \\
        SCIP(1.5s) & 200-400 & 1.3994 & $ 1.3952 \pm 0.0049 $ & $ 1.3636 \pm 0.0013 $ & $ 1.3608 \pm 0.0017 $ & $ 1.3629 \pm 0.0018 $ & \textbf{1.3536}  \\
        Gurobi(1.5s) & 200-400 & 1.4889 & $ 1.4879 \pm 0.0013 $ & $ 1.4754 \pm 0.0022 $ & $ 1.4721 \pm 0.0018 $ & $ 1.4697 \pm 0.0003 $  & \textbf{1.4684} \\
        \hline
    \end{tabular}%
    }
    \vspace{-15pt}
\end{table*}

\textbf{Solvers.} It has been shown that the greedy algorithm choosing a set that contains the largest weight of uncovered elements achieves an approximation ratio of $(1-\frac{1}{e})$~\citep{ratio}. Besides, the greedy algorithm has also been proved to be the best-possible polynomial-time approximation algorithm for maximum coverage unless P=NP~\citep{greedy}. So we choose it to serve as our first attack target. Besides, the maximum coverage problem is quite suitable to be formulated as an integer linear program (ILP). Using Google ORTools API~(\url{https://developers.google.com/optimization}), we transform Eq.~\ref{eq:mc} into the ILP form and attack three other general-purpose solvers: Gurobi (the SOTA commercial solver)~\citep{llc2020gurobi}, SCIP (the SOTA open-sourced solver)~\citep{scip}, and CBC (an open-sourced MILP solver written in C++)~\citep{cbc}, to find vulnerabilities in up-to-date solvers. These solvers are set the same appropriate time limit for a fair comparison.

\textbf{Calibrated time.} In line with~\cite{nair2020solving}, we use the calibrated time to measure the running time of task solving to reduce the interference of backend running scripts and the instability of the server itself. The main idea is to solve a small calibration MIP continuously and use the average running time to measure the speed of the machine when solving different instances. Details are given in Appendix \ref{app:Calibrated_time}.

\textbf{Attack model.} The MC problem can achieve a no-worse optimum if we add new elements into certain sets, as we can possibly cover more elements while not exceeding the threshold $k$. MC can be treated as a bipartite graph, where edges only exist between a set node and an element node inside the set. So our attack model is to add edges between set nodes and elements nodes, leading to a theoretically better optimum but can mislead the solvers. To reduce the action space, we only add edges for the unchosen sets, since adding edges for selected sets will not affect the solver's solution. An intuitive attack example is shown in Fig.~\ref{fig:MCexample} in Appendix \ref{app:examples}.

\textbf{Dataset.} For the MC problem, the distribution follows that in ORLIB~(\url{http://people.brunel.ac.uk/~mastjjb/jeb/info.html}). The dataset consists of three set-element pairs 100-200, 150-300, and 200-400, each with 50 training samples and 20 test samples.

\textbf{Results and analysis.} Tab.~\ref{tab:MCAttack} records the attack results of four target solvers. The RL based attack outperforms other methods across all cases in this scenario, showing the promising power of our feature extraction method and policy gradient model. Besides, the attackers have even found some instances for CBC where it cannot figure out feasible solutions under the given time limit, strongly demonstrating the vulnerability of the solvers. This effect also appears in the more challenging task Maximum Coverage with Separate Coverage Constraint MCSCC for Gurobi in Appendix \ref{app:experimentMCSCC}.

\section{Conclusion and Outlook}\value{-5pt}

We propose Robust Combinatorial Optimization (ROCO), the first general framework to evaluate the robustness of combinatorial solvers on graphs without requiring optimal solutions or the differentiable property. 
% The optimal solutions are hard to solve even in the median-size CO problems given acceptable time, and the differentiable property will exclude a large number of solvers such as heuristics, branch-and-bound solvers that have been widely deployed. 
Alleviating these two requirements makes it generalized and flexible, helping us conduct an extensive evaluation of the robustness of 14 unique combinations of different solvers and problems. Our experiment discovers problem instances that have no worse optimal costs but can significantly degenerate the solver's performance. The solution quality of the commercial solver Gurobi can degenerate by 20\% compared with the raw solution, under the same solving time, based on our discovered instances. 

The paradigm opens up new space for further research, including the aspects: 1) Utilizing the hard cases generated by our framework to design defense models, such as adversarial training and hyperparameter tuning. 2) Developing a more general attack framework for CO problems that cannot be represented by graphs. 3) Besides our studied tasks, the attacking idea of loosening constraints and lowering costs is also generalizable on a wide range of CO problems such as Maximum Cut and Minimum Vertex Cover, which need further evaluation on their robustness.
% 3) Exploring the running time attack to solve out optimal solutions using Gurobi, finding the key constraints that limit the solving time.

% \begin{ack}
% Use unnumbered first level headings for the acknowledgments. All acknowledgments
% More information about this disclosure can be found at: \url{https://neurips.cc/Conferences/2022/PaperInformation/FundingDisclosure}.

% Do {\bf not} include this section in the anonymized submission, only in the final paper. You can use the \texttt{ack} environment provided in the style file to autmoatically hide this section in the anonymized submission.
% \end{ack}

\bibliographystyle{abbrv}
\bibliography{acmart}
\clearpage

\appendix

\section{Traditional Attack Algorithm} \label{app:SA}

As an example of traditional attack algorithm, we list the pseudo code of SA in Algorithm~\ref{alg:SA}.

\begin{algorithm2e}
{\small{
    \caption{\textbf{Simulated Annealing (SA) Attack}}
    \label{alg:SA}}}
    \KwIn{Input problem $Q$; solver $f_\theta$; max number of actions $K$; action sample number $M$;\\
          \qquad \quad Temperature decay $\Delta T$; coefficient $\beta$; bias $eps$.}% bias $eps$
    $Q^0 \leftarrow Q$; $Q^* \leftarrow Q^0$;  $T \leftarrow 1$;  \algcomment{initial temperature }\\
    \For {$k \leftarrow 1..K$}
    {
        {\rm flag = False;} \algcomment{if action is available}\\ 
        \For {$i \leftarrow 1..M$}{
            Random sample an edge $(x,y)$ in edge candidates of $Q^{k-1}$;  \\
            $\tilde{Q} \leftarrow $ add/delete the edge $(x,y)$ in $Q^{k-1}$; \algcomment{new state by tentative action} \\ 
            $P = \exp(\frac{\beta (c(f_\theta,\tilde{Q}) - c(f_\theta,Q^{k-1}) + eps}{T})$; \algcomment{action acceptance probability}\\ 
            \If{$ \mathrm{Random} (0,1) \leq P$}{
                flag = True;  $Q^{k} \leftarrow \tilde{Q}$; $Q^* \leftarrow Q^k$\\
                \textbf{break};
            }
        }
        \If{{\rm flag = False}}{
            \textbf{break};
        }
        $T = T \cdot \Delta T$; 
    }
    \KwOut{Problem instance $\mathcal{G}^*$.}
\end{algorithm2e}

\section{Experiments-Task IV: Maximum Coverage with Separate Coverage Constraints} \label{app:experimentMCSCC}
We also study a more complicated version of the former MC problem: the maximum coverage problem with separate coverage constraint (MCSCC), which is NP-hard as proved in Appendix~\ref{app:MCSCCNPC}. Specifically, the elements $\mathcal{E}$ can be classified into the black ones $\mathcal{B}$ and the white ones $\mathcal{W}$ (i.e. $\mathcal{E}=\mathcal{W}\cup\mathcal{B},\mathcal{W}\cap\mathcal{B}=\emptyset$). And our goal is to \textbf{select a series of sets $\mathcal{S}' \subseteq \mathcal{S}$ to maximize the coverage of black element weights while covering no more than $k$ white elements.} The problem can be represented by a bipartite graph, where edges only exist between a set node and an element node covered by the set. The problem can be formulated as:
\begin{equation}\label{eq:MCSCCForm}
\begin{aligned}
%\max_R \sum_{b_j \in \mathcal{B}}v(b_j) \times \mathbb{I}(b_j \in \bigcup_{r_i \in R}B(r_i)) \qquad s.t. \quad |\bigcup_{r_i \in R}W(r_i)| \leq K 
\max_{\mathcal{S}'} \sum_{b \in \mathcal{B}}w(b) \times \mathbb{I}(b \in \bigcup_{S_i \in \mathcal{S}'}C^+(S_i))\\ 
s.t. \quad \mathcal{S}'\subseteq\mathcal{S},\ |\bigcup_{S_i \in \mathcal{S}'}C^-(S_i)| \leq k
\end{aligned}
\end{equation}
where $w(\cdot)$ denotes the weight of a certain element, $C(\cdot)$ denotes the set of elements covered by a set, and $C^+(\cdot)$ and $C^-(\cdot)$ denotes the set of elements in $C(\cdot)$ with black and white labels, respectively. %$B(\cdot)$ and $W(\cdot)$ represent the true label of events covered by a rule being black or white, which are covered by a rule that deems as black events.

%\subsubsection{Implementation Details}\vspace{-8pt}
\textbf{Solvers.} As an enhanced version of MC, the MCSCC problem is very challenging and here we propose three different solvers as the target for attacking. First, the trivial {Local} algorithm iterates over the sets sequentially, adding any sets that will not exceed the threshold $k$. Second, we adopt a more intelligent {Greedy Average} algorithm that always chooses the most cost-effective (the ratio of the increase of black element weights to the increase of a number of covered white elements) set at each step until the threshold $k$ is exceeded. Third, we formulate the problem into the standard ILP form (details in Appendix~\ref{app:MCSCCILP}) and solve it by Gurobi~\citep{llc2020gurobi}.

%\textbf{2) Attack model.} Suppose a criminal is trying to attack the banking system, one feasible and convenient way is to report some white events. Intuitively, when a white event is incorrectly classified  as a black event, the MCSCC problem will achieve no worse optimum $f^*(Q')$, since we can possibly cover more white events while not exceeding the threshold. Besides, we can shift our attention from the nodes to the edges. For the attacker, he can also report that an event is of a certain criminal nature. That is, to add additional black edges, which will lead to a theoretically better optimum, to fool our solvers. Furthermore, to reduce action space, we only choose rules that haven't been selected as our start node (or it will be useless since adding edges to the selected rules would not change the solvers' solution). Here we use the latter attack method and the former one's results are shown in the Appendix~\ref{app:table}.

\textbf{Attack model.} Our attack model chooses to add non-existing black edges that connect sets to black elements, which can lead to a theoretically better optimum since we can possibly cover more black elements while not exceeding the threshold $k$. Further, in order to reduce the action space, we only select the unchosen sets, otherwise it will be useless since adding edges for selected sets will not affect a solver's output solution. A successful attack example is shown in Fig.~\ref{fig:MCexample} in Appendix \ref{app:examples}.

% \textbf{4) Graph embedding.} For the RL attack method, different from DAG and TSP, MCSCC has a unique bipartite graph structure. Therefore, we resort to SAGEConv~\citep{hamilton2018inductive}, which can handle bipartite data, for graph feature extraction. As input, we classify the nodes into three classes (rules, black events and white events) and associate them with three dimension one-hot tensors. Besides, we add one more dimension for event nodes, which records their amounts. Eq.~\ref{GraphEmbed} becomes:
% \begin{equation}
% \begin{aligned}
% \mathbf{n}_e = \mbox{SAGEConv}_1(I_r,I_e)&,\  \mathbf{n}_r = \mbox{SAGEConv}_2(I_e,I_r) \\
% \mathbf{g}_e = \mbox{AttPool}_1(\mathbf{n}_e)&,\  \mathbf{g}_r = \mbox{AttPool}_2(\mathbf{n}_r)
% \end{aligned}
% \end{equation}
% % Then it is easy to generalize the new extraction features to fit in Eq.~\ref{Value} and~\ref{Prob} (details in Appendix).
% For heuristic attackers, we specify each action to be adding a black edge, and remaining details are just similar to what we've discussed in Section~\ref{Heuristic}.

%\subsubsection{Experiment Results}\vspace{-8pt}

\begin{figure}[tb!]
    \centering
    \includegraphics[width=0.7\linewidth]{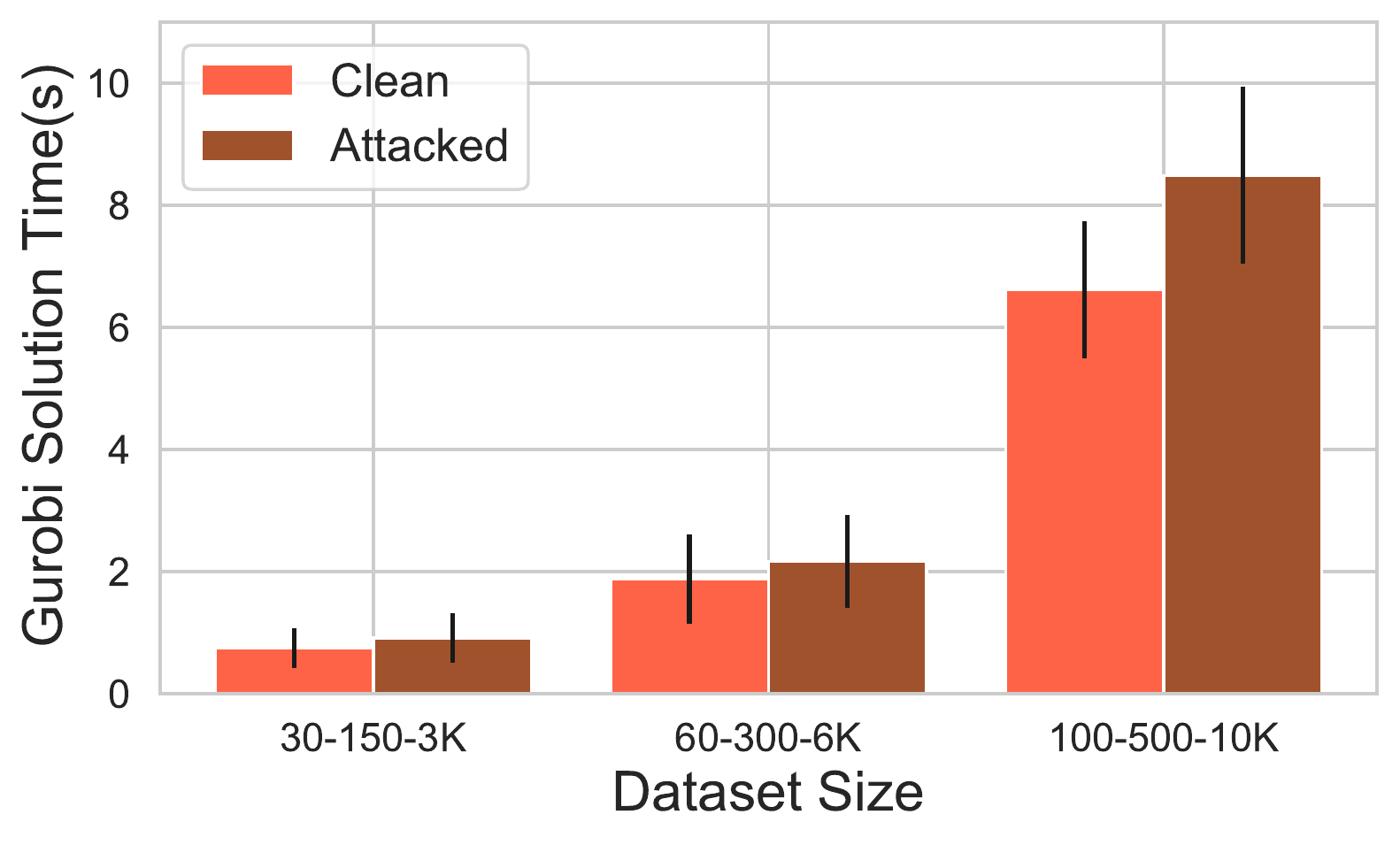}
    \vspace{-15pt}
    \caption{Gurobi's average time cost (no time limit) in solving clean or attacked MCSCC problems. Experiments are run on 3 datasets (20 instances) of different sizes.} 
    \label{fig:gurobi_time}
\end{figure}

\textbf{Dataset.} 
There is a lack of large-scale real-world datasets for MCSCC. Thus, we keep the data distribution similar to the MC problem and randomly generate the dataset~\cite{kwon2021matrix} for training. Specifically, the dataset consists of three configs of set-black element-white element: 30-150-3K, 60-300-6K and 100-500-10K, each with 50 training samples and 10 testing samples. More details about the dataset are presented in Appendix~\ref{app:MCSCC_Dataset}.
% We test on a real-world dataset from an Internet fintech company consisting of three rule and event pairs (i.e. 30-3K, 60-6K and 100-10K), each with 50 training samples and 20 testing samples. The distribution of transaction monetary values and the rule coverage are shown in Appendix~\ref{app:MCSCCData}. 

%Here time(s) refers to the time taken by Gurobi to find the optimal solution to MCSCC. The time change after attack and defense is given.
\begin{table*}[tb!]
    \centering
    \caption{Evaluation of solvers for MCSCC (notations are kept in line with MC in Table~\ref{tab:MCAttack}). Since MCSCC incorporates more constraints and is more challenging, our ROCO attackers demonstrate a more significant effect compared to MC, suggesting that solvers may be more vulnerable on CO problems that are inherently more challenging (e.g.\ with more complicated constraints). The greedy average algorithm is the most robust probably due to its worst-case guarantee known as submodularity~\cite{ratio}.}\label{tab:MCSCCAttack}
    \resizebox{\textwidth}{!}
    {
    \begin{tabular}{r|c|c|ccccc}
        \hline
        % Table generated by Excel2LaTeX from sheet 'MGM-willow'
        \multirow{2}[0]{*}{Solver} & Problem Size: & \multirow{2}[0]{*}{Clean Weight $\uparrow$} & \multicolumn{4}{c}{Attack Method (Attacked Weight)}  \\
         &  \#set-\#b elem.-\#w elem. & & Random & ROCO-RA & ROCO-OG & ROCO-SA & ROCO-RL   \\
        \hline
        Local Search & 30-150-3K  & 9.5713  & $ 9.5474 \pm 0.0325 $ & $ 9.4966 \pm 0.0057 $ & $ 9.4976 \pm 0.0105 $ & $ 9.4899 \pm 0.0029 $ & $ \textbf{9.4861}$   \\
        Greedy Average  & 30-150-3K  & 18.0038 & $ 18.2541 \pm 0.1710 $ & $ 17.5141 \pm 0.0288 $ & $ 17.4331 \pm 0.0414 $ & $ 17.5177 \pm 0.0396 $ &  $ \textbf{17.1414}$  \\
        Gurobi(1s) & 30-150-3K  & 18.8934 & $ 19.2297 \pm 0.1795 $ & $ 16.9266 \pm 0.2135 $ & $ 15.4132 \pm 0.3552 $ & $ 15.3055 \pm 0.3684 $ &  $\textbf{14.2740}$ \\
        \hline
        Local Search & 60-300-6K  & 24.9913 & $ 25.0413 \pm 0.0625 $ & $ 24.8738 \pm 0.0100 $ & $ \textbf{24.7914} \pm \textbf{0.035} $ & $ 24.8189 \pm 0.0375 $   &  $24.8014$  \\
        Greedy Average & 60-300-6K  & 43.1625 & $ 43.1539 \pm 0.0475 $ & $ 42.7697 \pm 0.0388 $ & $ 42.7611 \pm 0.0475 $ & $ 42.7222 \pm 0.0388 $ & $\textbf{42.1741}$ \\
        Gurobi(2s) & 60-300-6K  & 41.1828 & $ 41.8005 \pm 0.4407 $ & $ 38.2382 \pm 0.3459 $ & $ 37.3322 \pm 0.4201 $ & $ 38.2918 \pm 0.3665 $ & $\textbf{36.0432}$  \\
        \hline
        Local Search & 100-500-10K & 22.9359 & $ 22.9634 \pm 0.0482 $ & $ 22.7616 \pm 0.0206 $ & $ 22.6538 \pm 0.0183 $ & $ 22.7455 \pm 0.0206 $ &  $\textbf{22.5804}$  \\
        Greedy Average & 100-500-10K  & 51.3905 & $ 51.4008 \pm 0.0308 $ & $ 50.7481 \pm 0.0719 $ & $ \textbf{50.5169} \pm \textbf{0.1747} $ & $ 50.6865 \pm 0.0411 $  & $50.5631$  \\
        Gurobi(5s) & 100-500-10K  & 49.3296 & $ 52.8192 \pm 0.3803 $ & $ 48.1375 \pm 0.3597 $ & $ 47.4386 \pm 0.4933 $ & $ 49.2013 \pm 0.2467 $ &  $\textbf{43.4866}$ \\
        \hline
    \end{tabular}%
    }
\end{table*}

\textbf{Results and analysis.} Tab.~\ref{tab:MCSCCAttack} shows the attack results on our simulated dataset. Both traditional and RL approaches achieve significant attack effects, while RL outperforms the others in most cases (especially for {Gurobi}).
We can see from Tab.~\ref{tab:MCSCCAttack} that the degradation of Gurobi after attack is much more obvious than that in MC, possibly because the computation progress of MCSCC can include more vulnerable steps than MC. Fig.~\ref{fig:gurobi_time} shows the average time cost (no time limit) for Gurobi in solving MCSCC problems. The figure illustrates that solving time on the attacked instances is longer than those on the clean ones by a large margin, proving the promising power of our attack framework.

% \textbf{Results for attack and defense.} Table~\ref{tab:MCSCCDefense} records the results of attack and defense experiments on MCSCC problems. Experiments are conducted on the same test set. In general, the defender can compensate for the damage of attack effectively and obtain an even better solution than the baseline in some cases. Besides, as a commercial solver, {Gurobi} should be able to obtain optimal solutions if in sufficient time (assuming we have unlimited computational resources). So we record the time for {Gurobi} to find the optimal solution under attack and defense. The result is shown in Fig.~\ref{fig:gurobi}, where Gurobi's solution time after attack (defense) significantly increases (decreases). This inspires us to attack toward the solvers' solution time in future work.

\section{Graph Embedding for specific Tasks}
We discuss specific graph embedding metrics for different kinds of CO problems for reproducibility.
\label{app:graph_embedding}
\subsection{TASK I: DIRECTED ACYCLIC GRAPH SCHEDULING}
Since the task is a directed acyclic graph, we use {GCN} to encode the state in the original graph and its reverse graph with inversely directed edges separately. Then we concatenate the two node embeddings and use an attention pooling layer to extract the graph-level embedding for Eq.~\ref{GraphEmbed}:
\begin{equation}
    \mathbf{n} = [\mbox{GCN}_1(\mathcal{G}) || \mbox{GCN}_2(\mbox{reverse}(\mathcal{G}))],\ 
    \mathbf{g} = \mbox{AttPool}(\mathbf{n}).
\end{equation}

\subsection{Task II: Asymmetric Traveling Salesman Problem}
 Since the graph is fully connected, we use GCN to encode the state in the graph. Then we use an attention pooling layer to extract the graph-level embedding. Eq.~\ref{GraphEmbed} becomes: 
\begin{equation}
 \mathbf{n} = [\mbox{GCN}(\mathcal{G})],\ \mathbf{g} = \mbox{AttPool}(\mathbf{n}).    
\end{equation}

\subsection{Task III: Maximum Coverage}
For the RL attack method, different from DAG and ATSP, MC has a unique bipartite graph structure. Therefore, we resort to SAGEConv, which can handle bipartite data, for graph feature extraction. As input, we classify the nodes into two classes (subsets $I_s$ and elements $I_e$) and associate them with two dimension one-hot tensors. Besides, we add one more dimension for element nodes, which records their weights. Eq.~\ref{GraphEmbed} becomes:
\begin{equation}
\begin{aligned}
\mathbf{n}_e = \mbox{SAGEConv}_1(I_s,I_e)&,\  \mathbf{n}_s = \mbox{SAGEConv}_2(I_e,I_s) \\
\mathbf{g}_e = \mbox{AttPool}_1(\mathbf{n}_e)&,\  \mathbf{g}_s = \mbox{AttPool}_2(\mathbf{n}_s)
\end{aligned}
\end{equation}

\subsection{Task IV: Maximum Coverage with Separate Coverage Constraints}
The graph embedding mechanism for MCSCC is exactly the same as MC since they can both be represented by a bipartite graph except that we add one more one-hot dimension to distinguish between black and white elements.

\section{Successful Attack Examples} \label{app:examples}
\subsection{Task II: Asymmetric Traveling Salesman Problem}
\begin{figure}[h]
    \centering
    \includegraphics[width=0.9\linewidth]{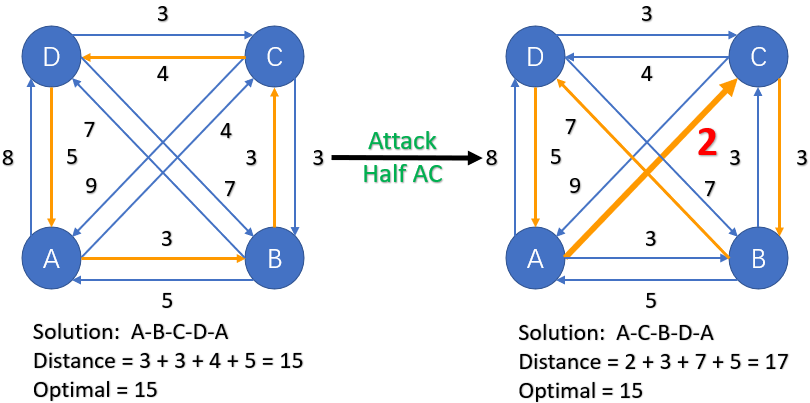}
    \caption{ATSP attack on {Nearest Neighbour} algorithm. The directed edge represents the distance between two nodes. The attack action on edge AC will cause 2 further distance.} 
    \label{fig:TSPexample}  
    \vspace{-10pt}
\end{figure}

A successful attack example is shown in Fig.~\ref{fig:TSPexample}.
\subsection{Task III: Maximum Coverage and Task IV: Maximum Coverage with Separate Coverage Constraints}
\begin{figure}[h]
    \centering
    \includegraphics[width=0.38\linewidth]{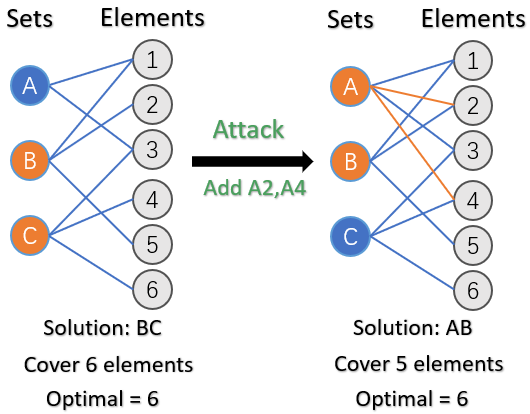}
    \includegraphics[width=0.61\linewidth]{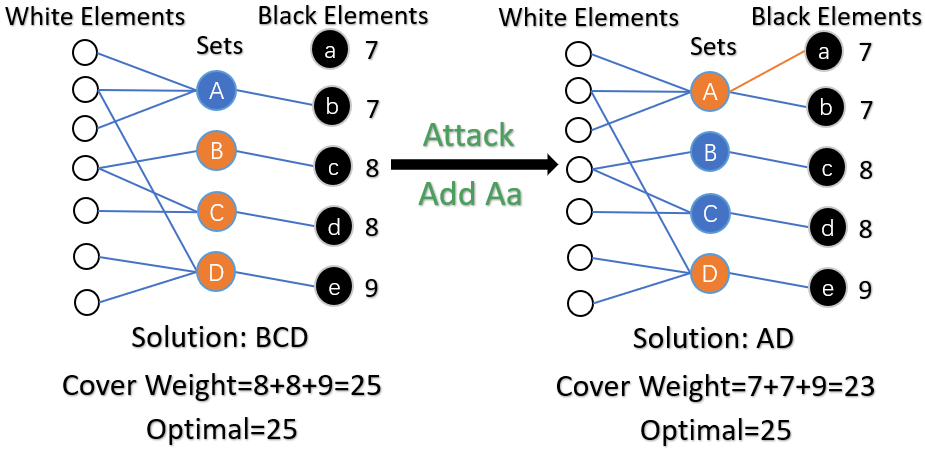}
    \caption{\textbf{Left}: MC attack on {Greedy} algorithm. The threshold $k$ of the sets is 2 and all the elements have unit weight. Attacks on edges A2 and A4 will cause a lower covered weight. \textbf{Right}: MCSCC attack on {Greedy Average} algorithm. The threshold $k$ of the white elements is $5$ and the value near the black element represents its weight. Attack on edge Aa causes a lower covered black element weight. \looseness=-1} 
    \label{fig:MCexample}
    \vspace{-10pt}
\end{figure}

Successful attack examples are shown in Fig.~\ref{fig:MCexample}.

\section{Formulas and Proofs} \label{app:prove}
\subsection{MCSCC NP-hard Provement}\label{app:MCSCCNPC}
We prove that the decision problem of \textbf{MCSCC} is NP-complete, thus the optimization problem of \textbf{MCSCC} in our paper is NP-hard. For the sake of our proof, here we redefine \textbf{MCSCC} and give definition of a traditional NP-Complete problem: \textbf{Set Covering}.

\textbf{The decision problem of MCSCC}. Given a series of sets $\mathcal{R}=\{R_1,R_2,...,R_n\}$ along with an element set consisted of white and black elements $\mathcal{E}=\mathcal{W}\cup\mathcal{B}$. The sets cover certain elements (white $C^-(R_i)$ or black $C^+(R_i)$) which have their weights $w(e)$. Does there exist a collection of these sets $\mathcal{R}'\subseteq\mathcal{R}$ to cover $\geq M$ black elements weight while influence no more than $k$ white elements?

\textbf{Set Covering}. Given a set $U$ of elements, a collection $S_1,S_2,\dots,S_n$ of subsets of $U$, and an integer $k$, does there exist a collection of $\leq k$ of these sets whose union is equal to $U$?

First, we need to show that \textbf{MCSCC} is NP. Given a collection of selected sets $\mathcal{R}'=\{R_1,R_2,\dots,R_m\}$, we can simply traverse the set, recording the covered black and white elements. Then, we can assert whether the covered black element number is no more than $k$ and fraudulent monetary value is no less than $m$.

Since the certification process can be done in $O(n^2)$, we can tell that \textbf{MCSCC} is NP. Then, for NP-hardness, we can get it reduced from \textbf{Set Covering}.

Suppose we have a \textbf{Set Covering} instance, we construct an equivalent \textbf{MCSCC} problem as follows: 
\begin{itemize}[leftmargin=*,itemsep=0pt,topsep=0pt]
    \item Create $|U|$ black elements with weight 1, corresponding to elements in \textbf{Set Covering}.
    \item Create $n$ sets, set $C^+(R_i)=S_i$.
    \item Connect each set to a different white element of weight 1.
    \item Set the white element threshold $k_w$ equal to the set number threshold $k_s$.
    \item Set black element weight target $M=|U|$.
\end{itemize} 

Suppose we find a collection of sets which meet the conditions of \textbf{MCSCC}, then we select subsets in \textbf{Set Covering} $S_i$ iff we select $R_i$. The total number of subsets is no more than $k_w(k_s)$ since the influenced white element number is equal to $|\mathcal{R}'|$. The subsets also cover $U$ since the covered black element weight (each with amount 1) is no less than $|U|$. Similarly, we can prove that we can find a suitable subset collection for \textbf{MCSCC} if we have found a collection of subsets that meet the conditions of \textbf{Set Covering}. So we can induce that $\textbf{Set Covering} \leq_p \textbf{MCSCC}$.

Thus, we have proved that the decision problem of \textbf{MCSCC} is NP-Complete.

\subsection{MCSCC ILP Formulation}\label{app:MCSCCILP}
As discussed in the main text, $\mathcal{B}$ denotes the set of black elements, $\mathcal{W}$ represents the set of white elements, $\mathcal{S}$ refers to the collection of sets and $\mathcal{E}$ denotes the element set. Using notations above, we can translate Eq.~\ref{eq:MCSCCForm} into standard ILP form as follows:
\begin{equation}
\begin{aligned}
\max \  \sum_{i \in \mathcal{B}} &Y[i] \times W[i], \quad s.t. \  \sum_{i \in \mathcal{W}} Y[i] \leq k \\
       \mbox{for}\  i=1\dots |\mathcal{E}|,\quad &(Y[i] - 0.5)(0.5 - \sum_{j=1}^{|\mathcal{R}|}X[j] \times \mathbb{I}(i\  \mbox{in}\  E[j])) \leq 0
\end{aligned}
\end{equation}
where $X[j] \in \{0,1\}$ denotes set $S_j$ is chosen (1) or not (0), $Y[i] \in \{0,1\}$ shows whether element $i$ has been covered by the chosen sets. Besides, $W[i] \in \mathbb{R}$ records the weight of the elements while $E[j] \subseteq \mathcal{E}$ is the corresponding elements of set $S_j$. The third equation ensures the element binary variable $Y[i]$ to be $1$ iff it has been covered by a certain set (if $\exists X[j]=1$ and element $i \in E[j]$, then the formula in the second bracket is negative, ensuring $Y[i]$ to be 1; else if element $i$ is not covered by any chosen sets, then the second formula is positive and $Y[i]$ must be 0).

\section{Calibrated Time} \label{app:Calibrated_time}
Considering the instability of the server and the interference of other running programs, calibrated time is designed to evaluate the server speed within a small period of  time. Specifically, we define the speed as the reciprocal of the Wall clock time to solve a given calibration MIP.
The calibration MIP is generated in the stable condition and solved $K$ times to record the basic speed $\text{speed}_{base}$. When evaluating the other instances, we will first evaluate the calibration MIP $K$ times to calculate the current speed $\text{speed}_{now}$. $K$ is the number of samples, which is set as $20$ in our experiments.

To fairly measure the solvers' performance under different computational resources, take MC problems as an example, a given calibration MIP is also a specified MC problem to design time limits for different server speeds. Specifically, the new time limit for the solver will be set as: 
\begin{equation}
    \text{New Time Limit} = \frac{\text{speed}_{base}}{\text{speed}_{now}} \times \text{Uniform Time Limit}
\end{equation}

\section{MCSCC Dataset} \label{app:MCSCC_Dataset}
Besides the distribution rules of MC, we generate the MCSCC dataset mainly by the following two rules: 1) The black element weights are uniformly distributed in the range $[0, 1]$. 2) The number of black elements is 5\% of white elements. 3) The covered black element weights and wight element numbers of sets follow a Gaussian distribution with a large standard deviation, to enlarge the differences between the sets.

\begin{figure}[tb!]
    \centering
    \includegraphics[width=0.9\linewidth]{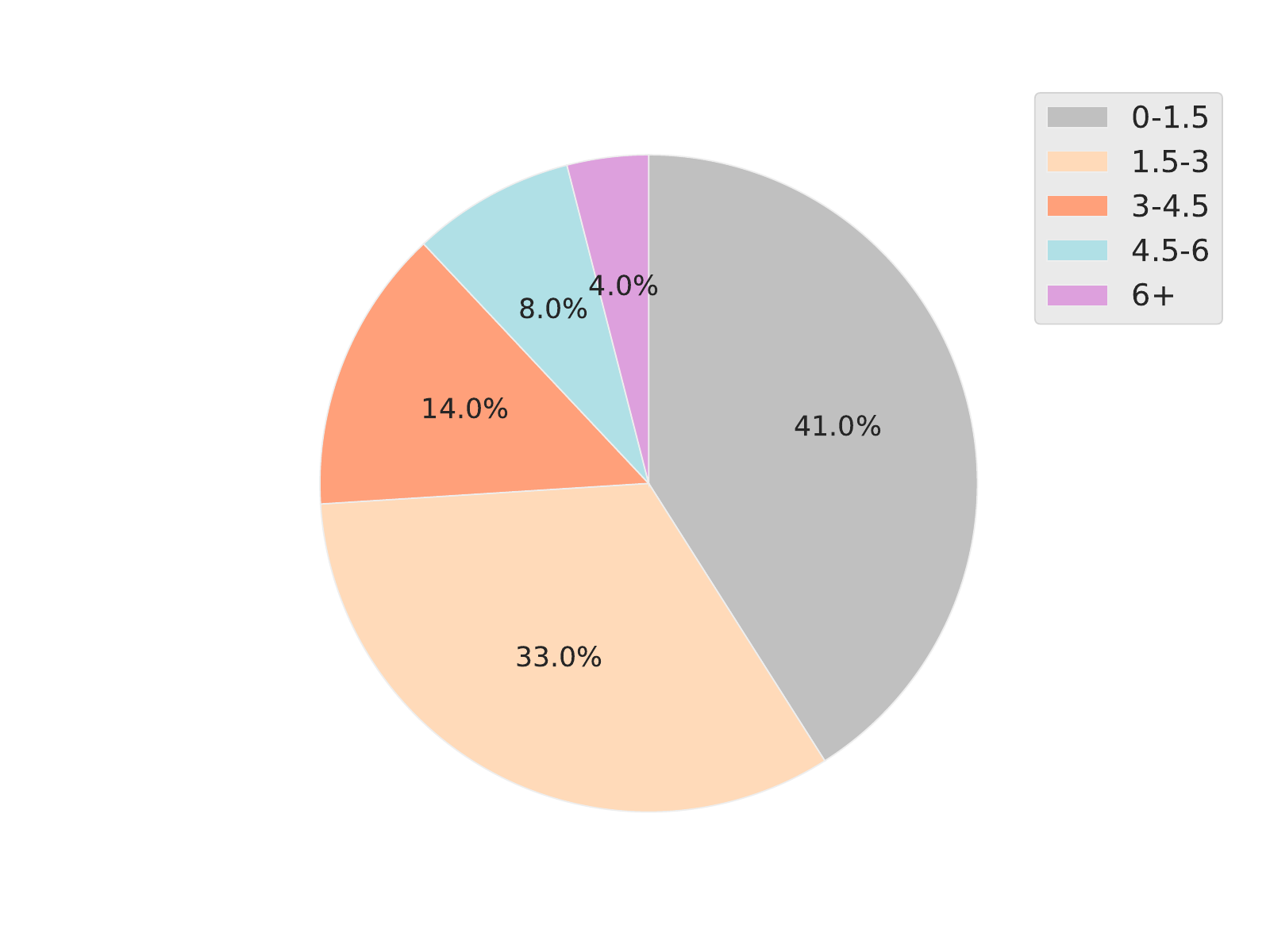}
    \vspace{-15pt}
    \caption{The distribution of black element weights covered by different sets. The $x-y$ pair in legends means the set can cover weights in range $[x,y]$.} 
    \label{fig:MCSCC_distrubition}
\end{figure}

Fig.~\ref{fig:MCSCC_distrubition} shows the black element weights covered by different sets in a 100-500-10K problem instance. As we can see, few well-designed sets can cover most weights while the others can be regarded as complementary to these sets. This design enlarges the difference between the different quality of solutions, therefore promoting the potential to be attacked.

% \begin{figure}[h]
% \vspace{-10pt}
% \centering
% \begin{minipage}[t]{0.49\textwidth}
% \centering
% \includegraphics[width=7.5cm]{iclr22/figures/Fraud_stragety.pdf}
% \caption{The distribution of fraudulent monetary value covered by different rules. The $x-y$ pair means the rules can detect amount in range $[x,y]$.} 
% \label{fig:fraud_strategy}  \vspace{-10pt}
% \end{minipage}
% \end{figure}

\section{Experiment Setups} 
\label{app:exp}
\textbf{Experiment environments.} DAG and ATSP experiments are run on a GeForce RTX 2080Ti while MC and MCSCC experiments are run on a GeForce RTX 3090 (20GB). We implement our models with Python 3.7, PyTorch 1.9.0 and PyTorch Geometric 1.7.2.

\textbf{RL settings.} Tab.~\ref{tab:RL_parameters} records the hyperparameters for RL during the training process. Trust region clip factor is a parameter in PPO agent to avoid model collapse. We also adopt some common policy-gradient training tricks like reward normalization and entropy regularization during training processes. \looseness=-1

\begin{table}[tb!]
    \centering
    \caption{RL parameter configuration in tasks DAG, ATSP, MC and MCSCC}
    \label{tab:RL_parameters}
   % \resizebox{0.9\columnwidth}{!}
    {
    \begin{tabular}{r|ccc}
        \hline
        % Table generated by Excel2LaTeX from sheet 'MGM-willow'
        Parameters & DAG & ATSP & MC \& MCSCC  \\
        \hline
        Actions\#  & 20 & 20 & 10  \\
        Reward discount factor & 0.95 & 0.95 & 0.95   \\
        Trust region clip factor & 0.1 & 0.1 & 0.1 \\
        GNN type  & GCN & GCN & SAGEConv \\
        GNN layers\#  & 5 & 3 & 3 \\
        Learning rate & 1e-4 & 1e-3 & 1e-3 \\
        Node feature dimensions\# & 64 & 20 & 16\\
        \hline
    \end{tabular}%
    }
\end{table}

\begin{table}[tb!]
	\centering
	\caption{Comparison of our attackers. Randomness means it will produce different results in different trials. Trained means whether the parameters are tuned by a training set.}
	\label{tab:timecomp}
	\centering
    \resizebox{0.6\columnwidth}{!}{
	 \begin{tabular}{r|ccc}
        \hline
        % Table generated by Excel2LaTeX from sheet 'MGM-willow'
       Attack Method & Randomness & Trained & Time Complexity \\
        \hline
        ROCO-RA & $\checkmark$ & & $O(NK)$\\ 
        ROCO-OG & $\checkmark$ & & $O(BMK)$\\
        ROCO-SA & $\checkmark$ & $\checkmark$ & $O(NMK)$\\
        ROCO-RL & & $\checkmark$ & $O(B^2K)$\\
        \hline
    \end{tabular}}
\end{table}

\textbf{Attacker hyperparameters.} For fair comparison of different attackers and the consideration of RL inference time, the hyperparameters are set to ensure similar evaluation time across different attack methods. According to the time complexity we calculate in Tab.~\ref{tab:timecomp}, here we specify the following parameters: number of iterations $N$, beam search size $B$ and number of different actions $M$ in each iteration. \looseness=-1

\newpage
\begin{itemize}[leftmargin=*,itemsep=0pt,topsep=0pt]
\item \textbf{DAG} : {ROCO-RA} $N = 30$; {ROCO-OG} $B = 3$, $M = 9$; {ROCO-SA} $N = 5$, $M = 6$; {ROCO-RL} $B = 3$.

\item \textbf{ATSP} : {ROCO-RA} $N = 130$; {ROCO-OG} $B = 5$, $M = 25$; {ROCO-SA} $N = 13$, $M = 10$; {ROCO-RL} $B = 5$.

\item \textbf{MC} : {ROCO-RA} $N = 220$; {ROCO-OG} $B = 6$, $M = 36$; {ROCO-SA} $N = 22$, $M = 10$; {ROCO-RL} $B = 6$.

\item \textbf{MCSCC} : {ROCO-RA} $N = 250$; {ROCO-OG} $B = 6$, $M = 36$; {ROCO-SA} $N = 25$, $M = 10$; {ROCO-RL} $B = 6$.
\end{itemize}

\section{Discussion about Limitations and Potential Negative Impacts}\label{app:discussion}

\paragraph{Limitations.} First of all, as our framework is designed for combinatorial optimization problems on graphs, it needs non-trivial reformulation to fit into CO problems cannot encoded by graphs. Secondly, the principle of designing our attack action is based on the graph modification to achieve no-worse optimum, which may not be applicable to all CO problems. 
% Some people may wonder whether all combinatorial optimization problems on the graph have the property. Fortunately, the actions to loosen the constrains or modify the objective function values tend to make sense although we can't exhaust all the questions. Apart from the four problems mentioned in our paper, the broadly studied CO problems on graphs Minimum Vertex Cover, Maximum Cut~\cite{khalil2017learning} also have the corresponding modifications to get the no-worse optimal solution. \looseness=-1

\paragraph{Potential Negative Impacts.}  As we propose the first general framework to evaluate the robustness of combinatorial solvers on graphs, this may be used by unscrupulous people to attack different solvers, increasing the burden on engineers to cope with the malicious attacks.

\end{document}